\documentclass[11pt,a4paper,leqno]{amsart}
\usepackage{amsmath,amsfonts,amssymb,amsthm}
\usepackage[english]{babel}
\usepackage{xcolor}
\usepackage{datetime}
\usepackage{enumerate}
\usepackage{caption}
\usepackage[colorlinks=true,linkcolor=blue, citecolor=blue, urlcolor=blue]{hyperref}%

\usepackage{pgfplots,tikz}
\usetikzlibrary{decorations.pathreplacing}
\usetikzlibrary{angles, arrows, patterns}
\usepgflibrary{patterns.meta}
\usepgfplotslibrary{fillbetween}

\usepackage[noadjust]{marginnote}
\usepackage{ifoddpage}

\newcommand{\mystar} 
{\checkoddpage\ifoddpage\normalmarginpar\else\reversemarginpar
\fi\marginnote{\n{$\star$}}}

\newcommand{\mystarn} 
{\checkoddpage\ifoddpage\reversemarginpar\else\normalmarginpar\fi\marginnote{\nr{\LARGE${
\Rrightarrow}$\quad}}}

\newcommand{\mystarb} 
{\checkoddpage\ifoddpage\reversemarginpar\else\normalmarginpar\fi\marginnote{\n{\LARGE${
\Rrightarrow}$\quad}}}

\newtheorem{thm}{Theorem}[section]

\newtheorem{cor}[thm]{Corollary}
\newtheorem{lemma}[thm]{Lemma}

\newtheorem{prop}[thm]{Proposition}
\numberwithin{equation}{section}

\theoremstyle{definition}
\newtheorem{remark}[thm]{Remark}

\def\R{\mathbb{R}}
\def\N{\mathbb{N}}
\def\Z{\mathbb{Z}}

\def\P{\mathbf{P}}
\def\E{\mathbf{E}}
\def\V{\mathbf{V}}

\newcommand{\uno}{\mathop{\textup{\large\textbf{1}}}}
\newcommand{\unotexto}{\mathop{\textbf{1}}}

\newcommand{\K}{\mathcal{K}}
\newcommand{\B}{\mathcal{B}}
\newcommand{\A}{\mathcal{A}}

\newcommand{\n}[1]{\textcolor{magenta}{#1}}

\pgfplotsset{compat=1.18}

\begin{document}

\title{Counting coprime pairs in random squares}

\author{Jos\'e L. Fern\'andez}
\address{Departamento de Matem\'aticas, Universidad Aut\'onoma  de Madrid, 28049 Madrid, Spain.}
\email{joseluis.fernandez@uam.es, \, jlfernandez@akusmatika.org}

\author{Pablo Fern\'andez}
\address{Departamento de Matem\'aticas, Universidad Aut\'onoma  de Madrid, 28049 Madrid, Spain.}
\email{pablo.fernandez@uam.es}

\thanks{Both authors are partially supported by Fundaci\'on Akusmatika.}

\begin{abstract}Extending the classical Dirichlet's density theorem on coprime pairs, in this paper we describe completely the probability distribution of the number of coprime pairs in random squares of fixed side length in the lattice $\N^2$. The  limit behaviour of this distribution as the side length of the random square tends to infinity is also considered. 
\end{abstract}




	\maketitle

\setcounter{tocdepth}{2}


\section{Introduction} The classical Dirichlet density theorem, coming all the way from \cite{Dirichlet}, claims that the proportion of coprime pairs in $\N_n\times\N_n$ converges to $6/\pi^2$ as $n \to \infty$:
\begin{equation}\label{eq:lim-dirichlet}
\lim_{n \to \infty}\,\frac{1}{n^2}\,\#\{(a,b) \in \N^2: 1 \le a,b \le n\, \mbox{and} \, \gcd(a,b)=1\}=\frac{6}{\pi^2}\cdot
\end{equation}

Informally (or formally), this theorem  claims that the probability that $\gcd(i,j)=1$, for a point~$(i,j)$ chosen at random and ``uniformly'' in $\N^2$, is $6/\pi^2$. ($\#A$ means ``number of elements of the set $A$''. Here and hereafter, for integer $n \ge 1$, we denote $\N_n:=\{1, \ldots, n\}$.)

The distribution of coprime pairs appear to be quite regular throughout the lattice  $\N^2$. This is reflected, for instance, in the fact that for a variety of random walks in $\N^2$, almost surely and asymptotically, the average time that the walker has coprime coordinates is again~$1/\zeta(2)$, see, for instance, \cite{CFF} for (regular) random walks and \cite{FF} for the so called P\'olya walks. See also Section 5 in \cite{FFracsam}.

Motivated by this anticipated  regularity, in this note, instead of drawing at random a point $(a,b) \in \N^2$  and checking whether it is a coprime pair or not, we select 
 at random a whole (square) \emph{window of fixed side length} $M \ge 1$, that is, a square of $M\times M$ lattice points, and count the number of coprime pairs within that  window. Our aim is to study the probability distribution  of this counting function. The case $M=1$ would correspond to Dirichlet's setting.

\begin{figure}[h]
\centering\resizebox{6.5cm}{!}{\begin{tikzpicture}[scale=1.0]

\draw[->,ultra thick] (-6,-5)--(13,-5) node[right]{\Huge $x$};
\draw[->,ultra thick] (-6,-5)--(-6,9) node[above]{\Huge $y$};

%
%



\draw[step=1.0,gray,thin] (-6,-5) grid (12,8);

\draw [ultra thick, decorate, decoration = {brace, mirror, amplitude=13pt}] (12.8,-0.2) --  (12.8,7.2) node[pos=0.5,right=16pt] 
     {\Huge $M$};


    \foreach \x in {5,6,...,12} {
        \foreach \y in {0,1,...,7} {
            \fill[color=black] (\x,\y) circle (0.2);
        }
        }

        \fill[color=black] (4.0,-1) circle (0.2) node[below, yshift=-15pt]{\Huge$(a,b)$};

          \end{tikzpicture}}
          \caption{The $M\times M$ window with base point $(a,b)$.}
         \label{fig:window} \end{figure}
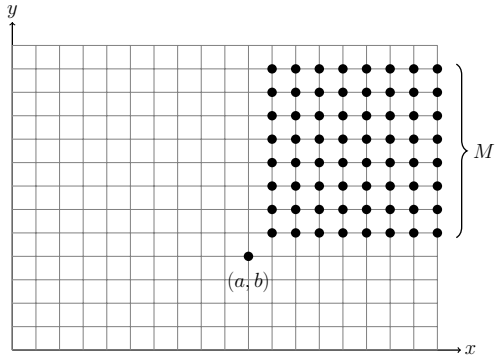

          Fix an integer $M\ge 1$, the side length of the window. For $(a,b)\in \N^2$, we let $Z_M(a,b)$ denote the number of coprime pairs $(i,j)$ within the square 
           of side~$M$ that `starts' from the point $(a,b)$, i.e., so that
\[
a+1\le  
 i\le a +  M \quad \mbox{and} \quad b+1\le  
  j\le b +  M.
\]
Thus, for $(a,b)\in \N^2$, we have that
\[Z_M(a,b)=\#\{(i,j): a+1\le  i\le  a+M, b+1\le j \le  b+M  \mbox{ and } \gcd(i,j)=1\}.\]
 Notice that the point $(a,b)$ is not exactly the south-west corner of the square. See Figure~\ref{fig:window}. The function $Z_M$ takes values in $\{0,\ldots, M^2\}$.

In \cite{SugitaTakanobu}, Sugita and Takanobu undertook the study of the function 
 $Z_M$, and proved that
the limit
\begin{equation}
\label{eq:lim-ST}
\lim_{n\to \infty} \frac{1}{n^2}\,\#\{(a,b) \in \N_n^2:   Z_M(a,b)=r\}
\end{equation}
exists for any $r$ such that $0\le r \le M^2$.

From a probabilistic point of view (see Section \ref{sec:probabilistic setting} for details), this means that the (random) variable $Z_M$ restricted to $\N_n^2$ converges in distribution,   as $n \to \infty$, to a (random) variable $Z_M^\star$ which takes values in $\{0, \ldots, M^2\}$.
The counting variable $Z_M^\star$ is a sum of $M^2$ Bernoulli variables, each of parameter $1/\zeta(2)$, which are not independent and actually have an interesting correlation structure, see Section~\ref{sec:correlation}.

(The existence of the limit distribution $Z_M^\star$ is also discussed in Theorem 1.1 of  the paper~\cite{Martineau} of Martineau.)

Sugita and Takanobu rather  considered the \textit{normalized} version $U_M$ of $Z_M$ given by
\[
U_M=M \Big(\frac{Z_M}{M^2}-\frac{6}{\pi^2}\Big),
\]
and  verified that the characteristic function of $U_M$ restricted to $\N_n^2$ converges pointwisely, i.e.,
\[
\lim_{n \to \infty} \frac{1}{n^2}\sum_{a,b \ge 1}^n e^{-i x \,U_M(a,b)} \quad \mbox{exists and it is finite for any $x \in \R$}.
\] Their very elegant argument is framed within the compactification of $\Z$ given by the ring of finite integral adeles.

\smallskip

In this paper, we show directly, 
 with tools from \emph{elementary number theory and basic probability theory}, that the limit \eqref{eq:lim-ST} exists for any $r$ such that $0 \le r \le M^2$, and provide \emph{formulas}
  for the probability distribution of the limit variable $Z_M^\star$ (see Theorem \ref{teor:distribution of ZM}).

We will also study the correlation structure of the summands of the counting function~$Z_M^\star$, to obtain estimates of its variance which imply 
 some limiting behaviour of $Z_M^\star$ as the side length  $M$ of the window tends to $\infty$: the variable $Z_M^\star/M^2$ tends in probability to the constant $1/\zeta(2)$ as $M \to \infty$, see Theorem~\ref{teor:limit Y_n/n^2}.

\smallskip

The \emph{contents of this paper} are organized as follows.
Section \ref{sec:background} contains some required background on probability and divisibility,
on arithmetic functions, and on the inclusion/exclusion principle. In Section \ref{sec:counting coprime pairs in windows},
we study the probability distribution of the variable~$Z_M$ and of its limit $Z_M^\star$,
while Section \ref{sec:correlation} is devoted to study the correlation structure of $Z_M^\star$. Finally, Section \ref{sec:limits} discusses the limiting behaviour of $Z_M^\star$ as the side length $M$ of the window tends to infinity and poses a few additional questions.

\section{Notation and some background}\label{sec:background}

\subsection{Some notations}

For integer $n \ge 1$, we write $\N_n$ to abbreviate  $\{1, \ldots, n\}$, and~$\N_n^2$ to shorten~$\N_n{\times}\N_n$. For $x>0$, $\lfloor x\rfloor$ and $\{x\}$ denote, respectively, the floor (integer part) of~$x$ and the fractional part of $x$, such that $\{x\}=x-\lfloor x\rfloor$.

For a pair of nonnegative integers $a$ and $b$, not both 0, the greatest common divisor $\gcd(a,b)$ of $a$ and $b$ is the largest nonnegative integer which divides both $a$ and $b$. Thus $\gcd(a,0)=a 
$, for integer $a \ge 1$. We shall find it convenient to follow the convention that $g(0,0)=0$.

For integers $d\ge 1$ and $a\in \Z$, we denote by $r_d(a)$ the remainder of dividing $a$ by $d$, or alternatively, the residue of $a$ mod $d$. We also denote by $I_d$ the function in~$\N$ given by $I_d(n)=1$ if $d \mid n$, and 0 otherwise, for each $n \in \N$.

The sets of all primes numbers is denoted by $\mathcal{P}$. The notations $\prod_{p}$, or $\prod_{p \ge M}$, etc., are short for $\prod_{p \in \mathcal{P}}$, $\prod_{p \in \mathcal{P},p\ge M}$, etc.


For a finite set $B$, we denote with $\#B$ (or with $|B|$) the number of elements of $B$.

If $\Omega$ is a certain reference set and $A$ is a subset of $\Omega$, by $\unotexto_A$ we denote the indicator function of $A$, i.e., the function in $\Omega$ such that $\unotexto_A(\omega)=1$, if $\omega \in A$, and $\unotexto_A(\omega)=0$, if $\omega \in \Omega \setminus A$. Observe that $\unotexto_{A\cap B}\equiv \unotexto_A \cdot \unotexto_B$.

In a generic probability space, we use $\P(A)$ to denote the probability of the event $A$, and $\E(X)$ and $\V(X)$ to denote expectation and variance of a 
 random variable $X$. Also, $\mathbf{cov}(X,Y)$ and $\rho(X,Y)$ will denote, respectively, 
  covariance and correlation coefficient of the random variables $X$ and $Y$.


Some additional notations for arithmetic functions will be introduced in Section~\ref{sec:arithmetic functions}.

\subsection{Probability and divisibility}\label{sec:probabilistic setting}

For each $n \ge 1$, we let $\P_n$ be the equidistributed probability in $\N_n^2$, thus, if $B \subset \N_n^2$, then $\P_n(B)=|B|/n^2$. For a set $B \subset \N^2$, we abbreviate and write $\P_n(B)=\P_n(B \cap \N_n^2)$. Also, $\E_n(X)$ and $\V_n(X)$ denote expectation and variance of a random variable $X$ defined in $\N_n^2$, and $\mathbf{cov}_{n}(X,Y)$ and $\rho_n(X,Y)$ denote covariance and correlation coefficient of the random variables $X$ and $Y$ defined in $\N_n^2$. We also abbreviate $\E_n(X)=\E_n(X\cdot \unotexto_{\N_n^2})$ and $\V_n(X)=\V_n(X\cdot \unotexto_{\N_n^2})$ for $X$ defined in the whole of $\N^2$, and analogously  for $\mathbf{cov}_{n}(X,Y)$ and $\rho_n(X,Y)$ for $X$ and $Y$ defined in the whole of $\N^2$.

Observe that for primes $p$ and $q$ (different or not),
\[
\#\{(a,b)\in \N^2_n: p \mid a \,\, \mbox{and}\,\, q \mid b\}=\Big\lfloor\frac{n}{p}\Big\rfloor \Big\lfloor\frac{n}{q}\Big\rfloor ,\]
and thus that
\[
\lim_{n \to \infty} \P_n\big(\{(a,b)\in \N^2: p \mid a \text{ and }  q\mid b\}\big)=\frac{1}{pq}\cdot
\]
In general, for any integers $u$ and $v$, we have that
\begin{equation}\label{eq:prop en pares de residuos}
\lim_{n \to \infty} \P_n\big(\{(a,b)\in \N^2: r_p(a)=u \text{ and }  r_q(b)=v\}\big)=\frac{1}{pq}\cdot
\end{equation}

The following elementary lemma describes the asymptotic independence in $\N^2$ of (joint) divisibility by different primes, which is general enough to cover the needs of 
 this paper.

\begin{lemma} \label{lemma:asympotic independence} Let $\{q_1, \ldots, q_R\}$ be a finite collection of \emph{distinct} primes. Consider integers $u_1, \ldots, u_R$ and $v_1, \ldots, v_R$ such that
$0\le u_j, v_j < 
 q_j$, for $j=1,\dots,R$ \textup{(}which play the role of collections of residues\textup{)}. 
Let
\[
\Gamma_j=\big\{(a,b)\in \N^2: r_{q_j}(a)=u_j\, \text{ and } \,r_{q_j}(b)=v_j\big\}, \quad  \text{for $j=1,\dots,R$.}
\]
Then
\[
\lim_{n \to \infty} \P_n\Big(\bigcap_{j=1}^R \Gamma_j\Big)=\prod_{j =  1}^R \frac{1}{q_j^2}=\prod_{j=1}^R \big(\lim_{n \to \infty} \P_n(\Gamma_j)\big).
\]
For $1\le S<R$, we have that
\[
\lim_{n\to \infty} \P_n\Big(\bigcap_{i=1}^S \Gamma_i\, \Big|\, \bigcap_{j=S + 1}^R \Gamma_j\Big)=\prod_{i=1}^S \frac{1}{q_i^2}\, \cdot
\]
 \end{lemma}

\begin{proof}
 The first statement follows directly from the Chinese reminder theorem and \eqref{eq:prop en pares de residuos}. The second, concerning conditional probability, follows since, from the first statement, we have that
 \[
\lim_{n \to \infty} \P_n\Big(\Big(\bigcap_{i=1}^S \Gamma_i \Big) \ {\textstyle\bigcap}\ \Big(\bigcap_{j=S + 1}^R \Gamma_j\Big)\Big)=\prod_{i=1}^S \frac{1}{q_i^2} \,\prod_{j=S + 1}^R \frac{1}{q_j^2}  ,
\]
and
\[
\lim_{n \to \infty} \P_n\Big(\bigcap_{j=S + 1}^R \Gamma_j\Big)=\prod_{j=S + 1}^R \frac{1}{q_j^2} \, \cdot\qedhere
\]
 \end{proof}

\subsection{Inclusion/exclusion arguments}

We shall resort a number of times to arguments of inclusion/exclusion type; we record next the classical 
 principle of inclusion/exclusion, and a few variants.

For $t\ge 1$, let $A_1, \ldots, A_t$ be subsets of a finite set $\Omega$.
Let $\P$ be a probability measure on the set $\Omega$ (over all subsets of $\Omega)$.
Then the inclusion/exclusion principle reads
\begin{equation*}
\Big|\bigcup_{j=1}^t A_j\Big|
=\sum_{1\le j \le t} |A_j| -\sum_{1\le i <j \le t} |A_i\cap A_j|+\cdots +(-1)^{t+1} |A_1\cap\cdots \cap A_t|
\end{equation*}
and 
\begin{equation}\label{eq:IEprob}
\P\Big(\bigcup_{j=1}^t A_j\Big)
=\sum_{1\le j \le t} \P(A_j) -\sum_{1\le i <j \le t} \P(A_i\cap A_j)+\cdots +(-1)^{t+1} \P(A_1\cap\cdots \cap A_t).
\end{equation}


We collect now convenient variations of the inclusion/exclusion principle:  the Schuette--Nesbitt formula, the Waring formula,  and the generating function approach,  which relate, for a given $r$,  the subset of $\Omega$ where a point lies in (exactly) $r$ of the subsets $A_j$ with the set of points lying in at least on $s$ of the~$A_j$, for each $s$. These variations are widely used, for instance,  in actuarial science. See 
  \cite{Gerber} and \cite[p. 89]{Gerberbook} for further details.

Let $C$ be the counting function of the $A_j$, i.e., \[C=\sum_{j=1}^t \uno\nolimits_{A_j}.\]
\begin{lemma}[Schuette--Nesbitt]\label{lemma:schuette-nesbitt} With the
 notations above,   for each integer $r$ such that $0\le r\le t$, we have that
\[
\begin{aligned}
\uno\nolimits_{\{C=r\}}(\omega)&=\sum_{s=r}^t (-1)^{s-r} \binom{s}{r} \bigg[\sum_{\substack{J\subset\{1,\ldots, t\},\\ |J|=s}} \,\,\prod_{j \in J} \uno\nolimits_{A_j}(\omega)\bigg]
\\
&=\sum_{s=r}^t (-1)^{s-r} \binom{s}{r} \bigg[\sum_{\substack{J\subset\{1,\ldots, t\},\\ |J|=s}} \,  \uno\nolimits_{\bigcap_{j \in J}A_j}(\omega)\bigg] , \quad \mbox{for each $\omega \in \Omega$}.
\end{aligned}
\]
%
\end{lemma}

Taking expectations,  Lemma \ref{lemma:schuette-nesbitt} gives the so-called Waring's formula:
\begin{equation}\label{eq:Waring1}
\P(C=r)=\sum_{s=r}^t (-1)^{s-r} \binom{s}{r} \sum_{\substack{J\subset\{1,\ldots, t\},\\ |J|=s}} \P\Big(\bigcap_{j \in J}A_j\Big),
\end{equation}
See, for instance p. 106 in Feller's book \cite{Feller}. 


For the counting variable $C=\sum_{1\le j \le t} \unotexto_{A_j}$, we have that
\[
\binom{C}{2}=\frac{C(C-1)}{2}=\sum_{1\le i<j\le t} \uno\nolimits_{A_i}\,\uno\nolimits_{A_j},
\]
 and, in general,
\[
\binom{C}{s}=\frac{C(C-1)\cdots (C-s+1)}{s!}=\sum_{1\le i_1<\cdots < i_s\le t} \uno\nolimits_{A_{i_1}}\cdots\uno\nolimits_{A_{i_s}} , \quad \mbox{for $1\le s \le t$}.
\]
(Notice that, for $s=0$,  we understand $\binom{C}{0}\equiv 1$ and also, consistently, that an empty intersection is the whole set: $\bigcap_{j \in \emptyset} A_j=\Omega$.) 
Thus,
\[
\E\Big(\binom{C}{s}\Big)=\sum_{\substack{J\subset \{1,\ldots, t\},\\ |J|=s}} \P\Big(\bigcap_{j\in J} A_{j}\Big) , \quad \mbox{for $0 \le s \le t$} ,
\]
and Waring's formula \eqref{eq:Waring1} can be rewritten as
\[
\P(C=r)=\sum_{s=r}^t (-1)^{s-r} \binom{s}{r} \,\E\Big(\binom{C}{s}\Big).
\]

\begin{lemma}[Inclusion/exclusion principle and probability generating functions]\label{lemma:inclusion/exclusion and pgf}
With the notations above,
\begin{equation}\label{eq:inclusion/exclusion and pgf}
\sum_{r=0}^t \P(C=r) \,z^r=\sum_{s=0}^t (z-1)^s \,\E\Big(\binom{C}{s}\Big) , \quad \mbox{for $|z| <1$}.
\end{equation}
\end{lemma}

\begin{proof} It follows from a direct change of order of summation:
\begin{align*}\sum_{r=0}^t \P(C=r) \,z^r&=\sum_{r=0}^t \sum_{s=r}^t (-1)^{s-r} \binom{s}{r}\,\E\Big(\binom{C}{s}\Big)\, z^r\\&=\sum_{s=0}^t \E\Big(\binom{C}{s}\Big) \sum_{r=0}^s (-1)^{s-r} \binom{s}{r} z^r=\sum_{s=0}^t (z-1)^s
\,\E\Big(\binom{C}{s}\Big).\qedhere
\end{align*}
\end{proof}

If the $A_j$ are \emph{exchangeable events} with respect to the probability $\P$, in the sense  that  for any $J \subset \{1, \ldots, t\}$ the probability $\P(\bigcap_{j \in J} A_j)$  depends only on $|J|$, and if we define $\alpha(s)=\P(\bigcap_{j \in J} A_j)$, for any $J \subset \{1, \ldots, t\}$ with $|J|=s$ and  $s \in \{0, \ldots, t\}$,  then~\eqref{eq:Waring1} and~\eqref{eq:inclusion/exclusion and pgf} reduce to
\[
\P(C=r)=\sum_{s=r}^t (-1)^{s-r} \binom{s}{r} \binom{t}{s} \alpha(s) , \quad \mbox{for $0\le r\le t$}.
\]
and
%
%
\[
\sum_{r=0}^t \P(C=r) \,z^r=\sum_{s=0}^t \binom{t}{s} \,(z-1)^s \,\alpha(s) , \quad \mbox{for $|z| <1$},
\]
respectively.

\subsection{Some results on arithmetic functions}\label{sec:arithmetic functions}
An arithmetic function $f\colon \mathbb{N}\to\mathbb{C}$ is termed \emph{multiplicative} if $f(1)=1$ and $f(n\cdot m)=f(n)\cdot f(m)$ for coprime $n$ and $m$, i.e., when $\gcd(n,m)=1$. If $f(n\cdot m)=f(n)\cdot f(m)$ holds for any $n,m\in\mathbb{N}$, then $f$ is termed \emph{completely multiplicative}. The multiplicative arithmetic function $f$ is called \emph{strongly multiplicative} if $f(p^a)=f(p)$ for all prime numbers $p$ and all natural numbers $a$.

For a multiplicative arithmetic function $f$ which is bounded (or simply, such that $|f(n)|=O_\varepsilon(n^\varepsilon)$, for every $\varepsilon >0$), its associated  Dirichlet series $L_f(s)$ admits an Euler product representation of the form
\[
L_f(s)=\sum_{n=1}^\infty \frac{f(n)}{n^s}=\prod_p \Big(1+\frac{f(p)}{p^s}+ \frac{f(p^2)}{p^{2s}}+\cdots\Big),\quad \mbox{for any $s \in \mathbb{C}$ such that $\Re s>1$} ,
\]
that reduces, in case $f$ is strongly multiplicative,  to
\begin{equation}\label{eq:Dirchlet series of strongly}
L_f(s)=\sum_{n=1}^\infty \frac{f(n)}{n^s}=\prod_p \Big(1+f(p)\, \frac{1}{p^s}\, \frac{1}{1-1/p^s}\Big), \quad \mbox{for any $s \in \mathbb{C}$ such that $\Re s>1$}.
\end{equation}

The \emph{$($Dirichlet$)$ convolution} $f\star g$ of two arithmetic functions $f$ and $g$ is given by
\[
(f\star g)(n)=\sum_{d\mid n} f(d)\, g(n/d)\quad\text{for $n\ge 1$.}
\]

For bounded arithmetic functions $f$ and $g$, their Dirichlet convolution can be bounded by $|(f\star g)(n)|\le C d(n)$, where $d(n)$ counts the number of divisors of the integer $n \ge 1$, and thus
$|(f\star g)(n)| =O_\varepsilon(n^\varepsilon)$, for any $\varepsilon>0$, according to Theorem 315 in \cite{HardyWright}.   The Dirichlet series of their convolution is the product of the individual Dirichlet series: $L_{f\star g}(s)=L_f(s)\cdot L_g(s)$, for $s \in \mathbb{C}$ such that $\Re s>1$.

The M\"{o}bius function $\mu$ is the multiplicative arithmetic function defined by  $\mu(1)=1$, and, for $n\ge 2$, by
\begin{equation}\label{eq:def de mu}
\mu(n)=
\begin{cases}
1,&\text{if $n$ is square-free and has an even number of prime factors,}
\\
-1,&\text{if $n$ is square-free and has an odd number of prime factors,}
\\
0,&\text{if $n$ is not square-free.}
\end{cases}
\end{equation}
Observe that $|\mu(n)|=1$, if $n$ is square free, and that $|\mu(n)|=0$, otherwise.

Recall that the Dirichlet series $L_\mu$ of the M\"{o}bius function $\mu$, and the Riemann zeta function, $\zeta(s)$, are related by
\[
L_\mu(s) =\sum_{n=1}^\infty \frac{\mu(n)}{n^s}=\prod_{p}\Big(1-\frac{1}{p^s}\Big)=\frac{1}{\zeta(s)} , \quad \mbox{for any $s \in \mathbb{C}$ such that $\Re s>1$}\,;
\]
in particular, we have that
\begin{equation}\label{eq:1/zeta(2) con Mobius}
L_\mu(2)=\sum_{n=1}^\infty \frac{\mu(n)}{n^2}=\prod_{p} \Big(1-\frac{1}{p^2}\Big)=\frac{1}{\zeta(2)}=\frac{6}{\pi^2}\approx 0\mbox{.}6079\dots
\end{equation}
The so called Feller--Tornier constant\footnote{In some instances, the Feller--Tornier constant is defined as $(1+\mathbf{F})/2$, with value $\approx 0\mbox{.}6613\dots$}, denoted here by $\mathbf{F}$, is given by
\begin{equation}\label{eq:Feller-Tournier}
\mathbf{F}=\prod_{p} \Big(1-\frac{2}{p^2}\Big)\approx 0\mbox{.}3226\dots
\end{equation}

\subsubsection{Inclusion/exclusion principle and M\"{o}bius function}\label{sec:inclusionexclusion mobius} We are going to encounter a few times the following situation: a finite set $\Omega$ endowed with a probability measure $\P$ defined for all subsets of $\Omega$,  a collection of subsets $A_p\subset \Omega$ indexed with prime numbers $p\in \mathcal{P}$, and a function $\beta$ so that $\P(A_p)=\beta(p)$ for any prime $p$, $\P(A_p\cap A_q)=\beta(pq)$ for any pair of distinct primes, and so on.

In this case, the inclusion/exclusion principle as in \eqref{eq:IEprob}, 
 combined with the codifying properties of the M\"{o}bius function, gives  us that
\[
\begin{aligned}
\P\Big(\Omega\setminus \bigcup_{p} A_p\Big)&=1-\P \Big(\bigcup_{p} A_p\Big)=1-\sum_p \P(A_p)+\sum_{p< q} \P(A_p\cap A_q)-\cdots
\\
&=1-\sum_p \beta(p)+\sum_{p< q} \beta(pq)-\cdots=
\sum_{\substack{h \ge 1, \\ \mbox{\tiny{square free}}}} \mu(h)\beta(h)=\sum_{h \ge 1} \mu(h) \beta(h),
\end{aligned}\]
which we register as
\begin{equation}\label{eq:inclusionexclusion mobius}
\P\Big(\Omega\setminus \bigcup_{p} A_p\Big)=\sum_{h \ge 1} \mu(h) \beta(h).
\end{equation}

\subsubsection{Ces\`{a}ro's identity}

The following standard identity, named after Ces\`{a}ro, see \cite{Cesaro}, is useful when dealing with sums over $\gcd$s.

\begin{lemma}[Ces\`{a}ro's identity]\label{lemma:cesaro identity} If $f$ is any arithmetic function, and for integers $A,B\ge 1$,
\[
\sum_{1\le i \le A,\,1\le j \le B} f(\gcd(i,j))=\sum_{k\ge 1} (f \star \mu)(k)\Big\lfloor \frac{A}{k}\Big\rfloor\Big\lfloor \frac{B}{k}\Big\rfloor .
\]
\end{lemma}

Actually, the sum on the right extends just up to $k \le \min\{A,B\}$.

\smallskip

If we apply Lemma \ref{lemma:cesaro identity} with  $A=B=n\ge 1$ and with the  function $f=\delta_1$, which is given by $\delta_1(k)=1$ if $k=1$, and 0 otherwise, then,  since $\delta_1\star \mu \equiv \mu$, we obtain that
\[
\#\{(i,j)\in \N_n^2: \gcd(i,j)=1\}=\sum_{k=1}^n \mu(k) \Big\lfloor \frac{n}{k}\Big\rfloor^2 .
\]
Using that $\lfloor n/k\rfloor=n/k-\{n/k\}$, a simple estimate shows that
\[
 \#\{(i,j)\in \N_n^2: \gcd(i,j)=1\}=n^2 \sum_{k=1}^\infty \frac{\mu(k)}{k^2}+O(n\ln n)\, \quad \mbox{as $n \to \infty$}.
\]
Therefore, by \eqref{eq:1/zeta(2) con Mobius}, we have that
\[
\frac{1}{n^2}\#\{(i,j)\in \N_n^2: \gcd(i,j)=1\}=\sum_{k=1}^\infty \frac{\mu(k)}{k^2}+O\Big(\frac{\ln n}{n}\Big)=\frac{1}{\zeta(2)}+O\Big(\frac{\ln n}{n}\Big),\quad \mbox{as $n \to \infty$,}
\]
and, thus, that
\begin{equation*}
\lim_{n \to \infty} \frac{1}{n^2}\,\#\{(i,j)\in \N_n^2: \gcd(i,j)=1\}=\frac{1}{\zeta(2)}\cdot
\end{equation*}
This observation is the Dirichlet density theorem anticipated in \eqref{eq:lim-dirichlet}.

\smallskip

We shall  need the following variation of Dirichlet's limit: 
\begin{equation}\label{eq:dirichlet general kl}
\lim_{n \to \infty} \frac{1}{n^2}\,\#\{(i,j)\in \N_n^2: \gcd(i + k,j + l)=1\}=\frac{1}{\zeta(2)} , \quad \mbox{for any  $(k,l)\in \N^2$}.
\end{equation}

To verify \eqref{eq:dirichlet general kl}, maintain $(k,l)$ fixed and let $B_n$ denote
\[
B_n=\{(i,j)\in \N^2: 1\le i \le k + n, 1\le j\le l + n\, \, \,\mbox{and} \,\, \,\gcd(i,j)=1\}.
\]
Now
\[
|B_n|=\sum_{m\ge 1} \mu(m) \Big\lfloor \frac{k + n}{m}\Big\rfloor\Big\lfloor \frac{l + n}{m}\Big\rfloor ,
\]
and arguing as above, we obtain that
\[
\lim_{n\to\infty} \frac{1}{n^2} \,|B_n|=\frac{1}{\zeta(2)}.
\]
Finally, 
 from the double inequality
\[
|B_n|-(k(l + n) + (k + n)l) \le \#\{(i,j)\in \N_n^2: \gcd(i + k,j + l)=1\}\le |B_n| ,
\]
we deduce  \eqref{eq:dirichlet general kl}.

\begin{proof}[Proof of Lemma \textup{\ref{lemma:cesaro identity}}] Fix first integers $C,D\ge 1$ and define, for every prime $p$, the set $A_p$ given by
 \[
 A_p=\{1 \le  i  \le  C, 1 \le j  \le  D: p\mid i \,\,\mbox{and}\,\, p \mid j\}.
 \]
 We have that $|A_p|=\lfloor C/p\rfloor \lfloor D/p\rfloor$, and that $|A_p\cap A_q|=\lfloor C/(pq)\rfloor \lfloor D/(pq)\rfloor$ for primes $p,q$, etc. From \eqref{eq:inclusionexclusion mobius}, with equidistributed probability, we obtain that
\[\#\{1 \le i \le  C, 1  \le  j  \le  D: \gcd(i,j)=1\}=CD-\Big|\bigcup_p A_p\Big|
=\sum_{h\ge 1} \mu(h)  \Big\lfloor \frac{C}{h}\Big\rfloor \Big\lfloor \frac{D}{h}\Big\rfloor.
\]
Now, using the above with $C=\lfloor A/m\rfloor$ and $D=\lfloor B/m\rfloor$, we deduce that
\begin{align*}
&\sum_{1\le i \le A,1\le j \le B}\! f(\gcd(i,j))=\sum_{m\ge 1} f(m) \, \, \#\{1\le i \le A, 1 \le j \le B: \gcd(i,j)=m\}
\\
&\qquad =\sum_{m\ge 1} f(m)\, \#\{1\le \tilde{i} \le A/m, 1 \le \tilde{j} \le B/m: \gcd(\tilde{i},\tilde{j})=1\}
\\
&\qquad =\sum_{m\ge 1} f(m) \sum_{h\ge 1} \mu(h)  \Big\lfloor \frac{\lfloor A/m\rfloor}{h}\Big\rfloor \Big\lfloor \frac{\lfloor B/m\rfloor}{h}\Big\rfloor
=\sum_{m\ge 1} f(m)  \sum_{h\ge 1} \mu(h)  \Big\lfloor \frac{A}{mh}\Big\rfloor \Big\lfloor \frac{B}{mh}\Big\rfloor
\\
&\qquad =\sum_{k\ge 1} \Big\lfloor \frac{A}{k}\Big\rfloor \Big\lfloor \frac{B}{k}\Big\rfloor \sum_{m,h\ge 1, mh=k} f(m)\,\mu(h)=\sum_{k \ge 1} (f \star \mu)(k)\Big\lfloor \frac{A}{k}\Big\rfloor \Big\lfloor \frac{B}{k}\Big\rfloor.
\end{align*}
We have used that if $n$, $m$ and $k$ are integers, then $\lfloor\lfloor n/m\rfloor/k\rfloor=\lfloor n/(mk)\rfloor$.) \end{proof}

\subsubsection{An auxiliary arithmetic function.}

In Section \ref{sec:correlation}, the function $\Upsilon $, which we are about to introduce, will be use to codify  the correlation structure of the counting function of coprime pairs in windows.

The arithmetic function $\Upsilon $ is defined  by
\begin{equation}\label{es:def de Upsilon}
\Upsilon (1)=1\quad \mbox{and} \quad \Upsilon (n)=\prod_{p  \mid  n} \frac{1-1/p^2}{1-2/p^2} , \quad \mbox{for $n \ge 2$}.
\end{equation}
The first few values of the function $\Upsilon (n)$ are
\begin{equation}\label{es:values of Upsilon}
1, \frac{3}{2},\frac{8}{7},\frac{3}{2},\frac{24}{23},\frac{12}{7},\frac{48}{47},\frac{3}{2},\frac{8}{7},\frac{36}{23},\frac{120}{119},\frac{12}{7},\frac{168}{167},\frac{72}{47},\frac{192}{161},\dots
\end{equation}


As all factors in the product defining $\Upsilon(n)$ are greater than 1, we have the bounds 
\begin{equation}\label{eq:cotas de Upsilon}
1 < \Upsilon (n)<\prod_{p} \frac{1-1/p^2}{1-2/p^2}=\frac{1}{\zeta(2) \cdot\mathbf{F}} \approx 1\mbox{.}88426\dots , \quad \mbox{for all $n  \ge 2$,}
\end{equation}
where $\mathbf{F}$ is the Feller--Tornier constant from \eqref{eq:Feller-Tournier}.



By convention, we extend $\Upsilon$ and define 
 $\Upsilon (0)=1/(\zeta(2)\cdot \mathbf{F})$


The function $\Upsilon $ is strongly multiplicative, due to its very definition \eqref{es:def de Upsilon}. Therefore,  by~\eqref{eq:Dirchlet series of strongly}, its associated Dirichlet series $L_{\Upsilon} $ 
 is given by
\[
L_{\Upsilon} (s)=\sum_{n=1}^\infty \frac{\Upsilon (n)}{n^s}=\prod_{p} \Big(1+\frac{1-1/p^2}{1-2/p^2}\,\frac{1}{p^s}\,\frac{1}{1-1/p^s}\Big) , \quad \mbox{for any $s \in \mathbb{C}$ such that $\Re s>1$}.
\]
In particular, for $s=2$, we have that
\begin{equation}\label{eq:value of Xi(2)}
L_{\Upsilon} (2)=\prod_{p} \Big(1+\frac{1-1/p^2}{1-2/p^2}\,\frac{1}{p^2}\,\frac{1}{1-1/p^2}\Big)=\frac{1}{\zeta(2)\, \mathbf{F}}\cdot
\end{equation}

\subsubsection{The convolution \texorpdfstring{$\,\Upsilon   \star\mu$}{\Upsilon\star\mu$}} We shall appeal to the convolution $\Upsilon  \star \mu$ of the function $\Upsilon  $ with the M\"{o}bius function $\mu$. Notice that the function $\Upsilon  \star \mu$ is multiplicative.
For prime $p$,
\[
(\Upsilon \star \mu)(p)=\sum_{d  \mid  p} \mu(d)\, \Upsilon(p/d)= \mu(1) \,\Upsilon (p)+\mu(p)\,\Upsilon (1)=\Upsilon (p)-1=\frac{1}{p^2 - 2},
\]
while for a prime power $p^k$ with~$k \ge 2$, 
\[
(\Upsilon \star \mu)(p^k)= \sum_{d  \mid  p^k} \mu(d)\, \Upsilon(p^k/d)   =\mu(1) \,\Upsilon (p^k)+\mu(p)\,\Upsilon (p^{k - 1})=\Upsilon (p) (\mu(1)+\mu(p))=0,
\]
using that $\Upsilon$ is strongly multiplicative. 
 Thus,
\begin{equation}\label{eq:def de upsilon*mu}
(\Upsilon \star\mu)(1)=1\quad \text{and}\quad (\Upsilon \star\mu)(n)=|\mu(n)| \prod_{p  \mid  n} \frac{1}{p^2 - 2} , \quad \mbox{for each $n \ge 2$}.
\end{equation}
Recall that $|\mu(n)|=1$ if $n$ is square free,  and is $0$  otherwise. The convolution $(\Upsilon \star\mu)$ is a non-negative 
 function, and in fact we have the bounds
\begin{equation}
\label{eq:bounds of OmegaStarMu}
|\mu(n)|\,\frac{1}{n^2}\le  (\Upsilon \star\mu)(n)\le  \frac{1}{\mathbf{F}}\, |\mu(n)|\,\frac{1}{n^2} , \quad \mbox{for any $n \ge 2$}.
\end{equation}
This follows from rewriting \eqref{eq:def de upsilon*mu} as
\[
(\Upsilon \star\mu)(n)=\frac{|\mu(n)|}{n^2}\, \prod_{p \mid  n} \frac{p^2}{p^2-2} , \quad \mbox{for any $n \ge 2$} ,
\]
and from observing that
\begin{equation*}
1< \prod_{p \mid  n} \frac{p^2}{p^2-2}<
\prod_{p} \frac{p^2}{p^2-2}=\frac{1}{\mathbf{F}} , \quad \mbox{for any $n \ge 2$}.
\end{equation*}
Notice that both inequalities in \eqref{eq:bounds of OmegaStarMu} are sharp: the constants 1 and $1/\mathbf{F}$ can not be improved; just take (big) prime numbers for the left inequality, and  primorials for the right inequality.



\subsubsection{Averages of\/ $\Upsilon $}

We consider now the average of $\Upsilon (\gcd(i,j))$ for $(i,j)\in \mathbb{N}_n^2$. 
Using Lemma \ref{lemma:cesaro identity}, we have that
\[
\sum_{1\le i,j \le n} \Upsilon (\gcd(i,j))=\sum_{k=1}^n (\Upsilon \star \mu)(k) \Big\lfloor \frac{n}{k}\Big\rfloor^2 , \quad \mbox{for any $n \ge 1$}.
\]
Writing $\lfloor n/k\rfloor=n/k+\{n/k\}$, and using the bound \eqref{eq:bounds of OmegaStarMu} and the values of $L_{\Upsilon} (2)$ from~\eqref{eq:value of Xi(2)} and $L_{\mu}(2)$ from~\eqref{eq:1/zeta(2) con Mobius}, we obtain that
\[
\begin{aligned}\sum_{1\le i,j \le n} \Upsilon (\gcd(i,j))&=n^2 \sum_{k=1}^\infty \frac{(\Upsilon \star\mu)(k)}{k^2}+O(n)
=n^2 \,L_{\Upsilon \star \mu}(2)+O(n)
\\
&=n^2 \,L_{\Upsilon} (2) \,L_\mu(2)+O(n)
=n^2 \,\frac{1}{\zeta(2)^2 \,\mathbf{F}}+O(n).
\end{aligned}
\]
Thus
\[
\E_n(\Upsilon  \circ \gcd)=\frac{1}{\zeta(2)^2 \,\mathbf{F}}+O\Big(\frac{1}{n}\Big) , \quad \mbox{ for  $n \ge 1$} ,
\]
and, in particular,
\[
\lim_{n \to \infty} \E_n(\Upsilon  \circ \gcd)=\frac{1}{\zeta(2)^2 \,\mathbf{F}}\, \cdot
\]

For 
 $\alpha, \beta \in (0,1]$, using that $\Upsilon$ is a positive function, we have that
 \[
 \sum_{\substack{1\le i \le \lfloor \alpha n\rfloor ,\\ 1\le j \le \lfloor \beta  n\rfloor }} \Upsilon (\gcd(i,j))
 \le
 \sum_{\substack{1\le i \le \alpha n,\\ 1\le j \le \beta  n}} \Upsilon (\gcd(i,j))
 \le
 \sum_{\substack{1\le i \le \lfloor \alpha n\rfloor +1,\\ 1\le j \le \lfloor \beta  n\rfloor +1}} \Upsilon (\gcd(i,j)),
 \]
and arguing as above, we deduce that
\[
\lim_{n\to \infty}\frac{1}{n^2}\sum_{\substack{1\le i \le \alpha n,\\ 1\le j \le \beta  n}} \Upsilon (\gcd(i,j)) = (\alpha \beta) \,\frac{1}{\zeta(2)^2 \,\mathbf{F}}  .
\]
More 
 generally, for any $\alpha, \beta,\gamma,\delta \in (0,1]$ such that $\alpha>\gamma$ and $\beta >\delta$,
\begin{equation}\label{eq:equi de upsilon}
\lim_{n\to \infty}\frac{1}{n^2}\sum_{\substack{\gamma n\le i \le \alpha n,\\ \delta n\le j \le \beta  n}} \Upsilon (\gcd(i,j))= (\alpha-\gamma)( \beta-\delta)  \,\frac{1}{\zeta(2)^2 \,\mathbf{F}}  \cdot
\end{equation}

The next lemma shows how well distributed the values $\Upsilon(\gcd(i,j))$ are. See \cite{FFequi1} for further connections between equidistribution and coprimality.

\begin{lemma}\label{lemma:promedio (1-x)(1-y) upsilon} If $f$ is a continuous function in the square $[0,1]^2$, then
\[
\lim_{n \to \infty} \frac{1}{n^2} \sum_{\substack{0\le i/n\le 1, \\ 0\le j/n\le 1}} f\Big(\frac{i}{n},\frac{j}{n}\Big)\, \Upsilon(\gcd(i,j))
=
\Big[\frac{1}{\zeta(2)^2\,\mathbf{F}}\Big]\, 
\int_0^1 \!\!\int_0^1
 f(x,y) \, dx \,dy.
\]
\end{lemma}

(The $i,j$ in the sum above are integers.)
This lemma claims that  the sequence of probability measures $\lambda_n$ in the square $[0,1]^2$ given by
\[
\lambda_n:=\frac{1}{\Lambda_n}\sum_{\substack{0\le i/n\le 1, \\ 0\le j/n\le 1}} \delta_{({i}/{n},{j}/{n})}\Upsilon(\gcd(i,j)) ,
\]
where $\delta_{(x,y)}$ denotes the point mass distribution at $(x,y)\in [0,1]^2$, and where
\[
\Lambda_n=\sum_{\substack{0\le i/n\le 1, \\ 0\le j/n\le 1}} \Upsilon(\gcd(i,j)),
\]
converges weakly, as $n \to \infty$, to the Lebesgue measure in the square $[0,1]^2$.


\begin{proof}[Proof of Lemma \textup{\ref{lemma:promedio (1-x)(1-y) upsilon}}]
Fix an integer $k\ge 1$. The $i,j$ in the sums below are always integers. Denote,  for integer $n > k$, 
\[
A_n:=\frac{1}{n^2} \sum_{\substack{0\le i/n\le 1, \\ 0\le j/n\le 1}} f\Big(\frac{i}{n},\frac{j}{n}\Big)\,\Upsilon(\gcd(i,j)).
\]

For integers $u$ and $v$ such that $0\le u,v <k$, denote
\[
A_n(u,v)=\frac{1}{n^2} \sum_{\substack{u/k\le i/n\le (u+1)/k, \\ v/k\le j/n\le (v+1)/k}} f\Big(\frac{i}{n},\frac{j}{n}\Big)\,\Upsilon(\gcd(i,j)).
\]
Let
\[
\phi(u,v)=\max\Big\{f(x,y):\frac{u}{k}\le x\le \frac{u + 1}{k} \mbox{ and }  \frac{v}{k}\le y\le \frac{v + 1}{k}\Big\} , \quad \mbox{for $0\le u,v <k$,}
\]
and
\[
B_n(u,v)=\phi(u,v)\, \frac{1}{n^2}\sum_{\substack{u/k\le i/n\le (u+1)/k, \\ v/k\le j/n\le (v+1)/k}} \Upsilon(\gcd(i,j)).
\]

Because $\Upsilon$ is a positive function, we have that
$A_n(u,v)\le B_n(u,v)$ for each $u$ and $v$. On account of \eqref{eq:equi de upsilon}, we have that
\[
\lim_{n \to \infty} B_n(u,v)=\phi(u,v)
\, \frac{1}{k^2 \,\zeta(2)^2 \,\mathbf{F}}\cdot
\]

Since
\[
A_n\le \sum_{\substack{0\le u <k,\\ 0\le v <k}}A_n(u,v)\le \sum_{\substack{0\le u <k,\\ 0\le v <k}}B_n(u,v) ,
\]
we have that
\[
\limsup_{n\to \infty} A_n\le \frac{1}{\zeta(2)^2 \mathbf{F}} \, \frac{1}{k^2} \sum_{\substack{0\le u <k,\\ 0\le v <k}}\phi(u,v).
\]

This last inequality is valid for any integer $k \ge 1$, and so from
\[\lim_{k \to \infty}\frac{1}{k^2}\sum_{\substack{0\le u <k,\\ 0\le v <k}}\phi(u,v)=
\int_0^1 \!\!\int_0^1
f(x,y) \,dx \, dy ,\]
(by definition of the Riemann integral), we deduce that
\[\limsup_{n\to \infty} A_n\le \frac{1}{\zeta(2)^2\,\mathbf{F}}
 \int_0^1 \!\!\int_0^1
f(x,y) \,dx \,dy.\]

Analogously, one obtains that
\[\liminf_{n\to \infty} A_n\ge \frac{1}{\zeta(2)^2\,\mathbf{F}}
  \int_0^1 \!\!\int_0^1
  f(x,y) \,dx \,dy.\qedhere\]
\end{proof}

The following particular example of the lemma above, with $f(x,y)=(1-x)(1-y)$, will be used later in this paper:
\begin{equation}\label{eq:promedio (1-x)(1-y) upsilon}
\lim_{n\to \infty} \frac{1}{n^2} \sum_{\substack{0\le i/n\le 1, \\ 0\le j/n\le 1}} \Big(1-\frac{i}{n}\Big)\Big( 1-\frac{j}{n}\Big)\,\Upsilon(\gcd(i,j))
=\frac{1}{4}\,\frac{1}{\zeta(2)^2\,\mathbf{F}} \cdot
\end{equation}

%
%

\section{Counting coprime pairs in windows}\label{sec:counting coprime pairs in windows}

In this section, we study the distribution of the variable $Z_M$ which counts coprime pairs in square windows of fixed side length $M$ to obtain Theorem  \ref{teor:distribution of ZM}, the main result of this paper.

\subsection{Indicator of coprime pairs}

We denote by $\mathcal{V}$ the set of points $(n,m)$ in $\N^2$ such that $\gcd(n,m)=1$, that is, so that $(n,m)$ is a coprime pair. Points of $\mathcal{V}$ are frequently called \emph{visible points} (from the origin $(0,0)$), see, for instance, \cite{HerzogStewart}.

We denote the indicator function of $\mathcal{V}$ by $V$, so $V(n,m)=1$ if $\gcd(n,m)=1$, and $V(n,m)=0$ otherwise.
We can write the function 
 $V$ as
\begin{equation}\label{eq:formula V(a,b)}
V(a,b)=\prod_{p } \big(1-I_p(a) \,I_p(b)\big) \quad\text{for $a,b\ge 1$}.
\end{equation}
Recall that $I_p(n)=1$ if the prime $p$ divides $n$, and $I_p(n)=0$ otherwise. For each $(a,b)\in \N^2$, all but a finite number of factors in the above expression are equal to 1. In fact, $V(a,b)=1$ if and only all factors are 1, or equivalently, if no prime $p$ divides both $a$ and $b$.

In probabilistic terms, and using the notation of Section~\ref{sec:probabilistic setting},
equation \eqref{eq:lim-dirichlet}, Dirichlet's density theorem,  translates into
\begin{equation}\label{eq:lim-dirichlet-prob}
\lim_{n \to \infty} \E_n(V)=\lim_{n\to \infty} \P_n(\mathcal{V})=\frac{1}{\zeta(2)}\cdot
\end{equation}

\subsection{Counting coprime pairs in windows}

Fix a side length $M \ge 1$, and denote by~$\K_M$ the square $\K_M=\{1, \ldots, M\}^2$ in $\N^2$.

For each point $(a,b)\in \N^2$, the window $\mathcal{W}_M(a,b)$ in the lattice $\N^2$ is the translation by~$(a,b)$ of the square~$\K_M$, that is,
\[
\mathcal{W}_M(a,b)=\{a +  1, \ldots, a +  M\}\times \{b +  1, \ldots, b +  M\}=(a,b) +  \K_M.
\]
See again Figure \ref{fig:window}.

We denote with $Z_M$ the function defined in $\N^2$ which at each $(a,b)$ gives the number of coprime pairs within the window $\mathcal{W}_M(a,b)$, i.e.,
\begin{equation}\label{eq:def of ZM-cero}
Z_M(a,b)=\sum_{(i,j)\in \mathcal{W}_M(a,b)} V(i,j).
\end{equation}
We prefer to write $Z_M$ in the form
\begin{equation}\label{eq:def of ZM}Z_M(a,b)=\sum_{(k,l)\in \mathcal{K}_M} V(a + k,b + l) ,\end{equation}
so as to display the function $Z_M$  as the sum of the $M^2$ functions
\[(a,b)\in \N^2 \mapsto V(a+k,b+l) ,\]
with $(k,l)$ running over the square $\mathcal{K}_M$. In principle, we have that $0 \le Z_M(a,b)\le M^2$, although the maximum value of $Z_M$ is usually quite smaller than $M^2$; see Section~\ref{section:pob distr of Z and Z*}.

\begin{remark}\label{remark:caso Z=0}
The case $Z_M(a,b)=0$ would correspond to a $M$-window which  is \emph{fully invisible} (from the origin), in the sense that all points of the window have coordinates that are not coprime. There are explicit constructions of arbitrarily large `invisible' squares; see, for instance, Theorem 5.29 in Apostol's book \cite{Apostol}, or \cite{HerzogStewart}. See also Section 5 in \cite{FFracsam} for some related questions.
\end{remark}

%
%
%
%

The following explicit formula for $Z_M(a,b)$ can be obtained by means of the inclusion/exclusion principle. Compare with Theorem 5 in~\cite{SugitaTakanobu}.

\begin{lemma}[A formula for $Z_M(a,b)$]\label{lemma:formula ZM} 
For $M\ge 1$ and $(a,b)\in\mathbb{N}^2$,
\begin{equation}\label{eq:formula ZM}
Z_M(a,b)=\sum_{d \ge 1} \mu(d)\Big\lfloor \frac{M+r_d(a)}{d}\Big\rfloor \Big\lfloor \frac{M+r_d(b)}{d}\Big\rfloor.
\end{equation}
\end{lemma}

\begin{proof}

Fix $M\ge 1$. For $d \ge 1$,  let
\[
C^{(d)}_M(a,b)=\{(i,j)\in (a,b)+\K_M: d\mid i\, \mbox{and} \, d \mid j\}.
\]
Now,
\[
Z_M(a,b)=M^2-\Big|\bigcup_{p\in \mathcal{P}} C^{(p)}_M(a,b)\Big|.
\]
The union above is in fact a finite union, since for $p> a+M$, the set  $C^{(p)}_M(a,b)$ is empty. Observe that, for prime $p$,
\[
|C^{(p)}_M(a,b)|=\Big\lfloor \frac{M+r_p(a)}{p}\Big\rfloor \Big\lfloor \frac{M+r_p(b)}{p}\Big\rfloor ,
\]
and that, for primes $p$ and $q$,
\[
\big|C^{(p)}_M(a,b)\cap C^{(q)}_M(a,b)\big|=|C^{(pq)}_M(a,b)|=\Big\lfloor \frac{M+r_{pq}(a)}{pq}\Big\rfloor \Big\lfloor \frac{M+r_{pq}(b)}{pq}\Big\rfloor ,
\]
and so on. Thus,  \eqref{eq:formula ZM} follows by \eqref{eq:inclusionexclusion mobius} of Section \ref{sec:inclusionexclusion mobius}.
\end{proof}

Going back to the probabilistic setting in $\N_n^2$, the random variable $Z_M$ is a sum of $M^2$ Bernoulli variables, see \eqref{eq:def of ZM}. Each one of them has, asymptotically as $n \to \infty$, parameter $1/\zeta(2)$, by
Dirichlet's result \eqref{eq:lim-dirichlet-prob}.
 All this readily gives, for the mean of $Z_M$, that
\begin{equation}
\label{eq:lim mean of ZM in Nn}
\lim_{n \to \infty} \E_n(Z_M)=\frac{M^2}{\zeta(2)}\cdot
\end{equation}
  But for the \emph{actual distribution} of $Z_M$, observe that those Bernoulli variables are not independent; in fact, they exhibit an interesting correlation structure, see Section \ref{sec:correlation}.

\smallskip
Our main interest here is to show that for each $M\ge 1$, the random variable $Z_M$ in $(\N_n^2, \P_n)$ converges in distribution, as $n\to\infty$, to a random variable $Z_M^\star$ taking values in $\{0,1,\dots,M^2\}$, i.e., to prove that the limit
\[
\lim_{n \to \infty} \P_n(Z_M=r)\,
=:
\P(Z_M^\star=r)
\]
exists for each $r$ such that $0 \le r \le M^2$  and that $\sum_{r=0}^{M^2} \P(Z_M^\star=r)=1$.



\smallskip

The analysis in the following Sections~\ref{sec:splitting} and~\ref{sec:distr of Z and Z*} will provide, see Theorem~\ref{teor:distribution of ZM}, an \emph{explicit} and \emph{computable}
formula for $\P(Z_M^\star=r)$ for each $M \ge 1$ and all $0 \le r \le M^2$.

\subsection{Splitting the square \texorpdfstring{$\K_M$}{K\_M}}\label{sec:splitting}

There is a natural interaction between the side length~$M$ of the square~$\K_M$ and divisibility properties of the points within the window $\mathcal{W}_M(a,b)$. The (simple) reason is that a prime $p\ge M$ cannot divide simultaneously $a+k$ and $a+k'$ if $k,k'\in\{1,\dots,M\}$. So, it will be most convenient to separate, for each $(a,b)\in \N^2$, the points of $\K_M$ into two classes, as follows.

Fix $M\ge 1$ and define
\[
P_M:=\prod_{p<M} p,
\]
 with $P_1=P_2=1$. We shall assign, to each $(a,b)\in \N^2$, a pair of complementary subsets, $\B_M(a,b)$ and $\A_M(a,b)$, within  the square~$\K_M$.

We denote by $\mathcal{B}_M(a,b)$ the subset of $\K_M$ which consists of those pairs $(k,l)\in\K_M$ such that both $a +  k$ and $b +  l$ are divisible by some prime $p <M$:
\begin{equation}\label{eq:def B_m}
\mathcal{B}_M(a,b)=
\bigcup_{p < M}\big\{(k,l) \in \K_M: a +  k \equiv 0 \!\!\mod p \ \text{  and } \ b +  l \equiv 0 \!\!\mod p\big\}.
\end{equation}
The set $\A_M(a,b)$ is just the complement of $\B_M(a,b)$ in $\K_M$
\begin{equation}\label{eq:def A_m}
\A_M(a,b):=\K_M\setminus\B_M(a,b).
\end{equation} 

Notice that $\mathcal{B}_1(a,b)=\mathcal{B}_2(a,b)=\emptyset$, and also that, accordingly, $\mathcal{A}_1(a,b)=\K_1$ and $\mathcal{A}_2(a,b)=\K_2$.

It is a relevant fact that the sets $
\mathcal{B}_M(a,b)$ and $\A_M(a,b)$ depend only on the \emph{collection} of (pairs of) residues $\{(r_q(a), r_q(b)): q \in P_M\}$, and thus, because of the Chinese remainder theorem, ultimately depend  only on the pair of residues \begin{equation}\label{eq:def of RM}
R_M(a,b):=(r_{P_M}(a),r_{P_M}(b)),
\end{equation}
Thus, there are~$P_M^2$ different possibilities for  the sets $\A_M(a,b)$ (and for the corresponding sets $\B_M(a,b)$).

We denote by $\Phi_M(a,b)$ the size of $\A_M(a,b)$:
\[
\Phi_M(a,b)=|\mathcal{A}_M(a,b)|.
\]
Arguing with the inclusion/exclusion principle as in \eqref{eq:inclusionexclusion mobius}, we obtain that
\begin{equation}
\label{eq:formula of PhiM(a,b)}
\Phi_M(a,b)=\sum_{1 \le d \,\mid P_M} \mu(d) \Big\lfloor \frac{M +r_d(a)}{d}\Big\rfloor \Big\lfloor\frac{M +r_d(b)}{d}\Big\rfloor.
\end{equation}
The sum above extends to the divisors $d$ of $P_M$.

It is always the case, although far from sharp,  that $\Phi_M(a,b)\le M^2 -\lfloor(M - 1)/2\rfloor^2$.

As mentioned before, 
 $\Phi_1(a,b)=1$ and $\Phi_2(a,b)=4$ for all $(a,b)\in\mathbb{N}^2$.

 For $M=3$,  we have that  $P_3=2$ and that
\begin{itemize}
\item $\Phi_3(a,b)=8$ if both $a$ and $b$ are even, i.e., if $R_3(a,b)=(0,0)$;
\item $\Phi_3(a,b)=5$ if both are odd, i.e.,  if $R_3(a,b)=(1,1)$;
\item $\Phi_3(a,b)=7$ if one is even and the other is odd, i.e., if $R_3(a,b)$ is $(1,0)$ or $(0,1)$.
\end{itemize}

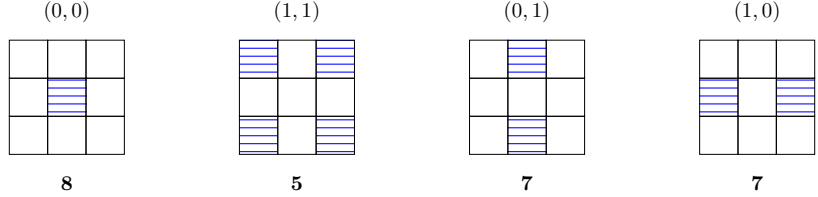
\begin{figure}[h]
\centering
\resizebox{0.7\linewidth}{!}{\begin{tikzpicture}[scale= 0.7]


              \fill [pattern= horizontal lines, pattern color =blue]  (1,1) rectangle (2,2);

            \draw (1.5, 3.75) node {$(0,0)$};
            \draw (1.5, -0.75) node {$\mathbf{8}$};

                \foreach \i in {0,1,2}{
                 \foreach \j in {0,1,2}{

              \draw  (\i,\j) rectangle (\i+1,\j+1);
              }
              }


  			 \foreach \i in {6,8}{
                 \foreach \j in {0,2}{

              \fill [pattern= horizontal lines, pattern color =blue]  (\i,\j) rectangle (\i+1,\j+1);
              }
              }

           \draw (7.5, 3.75) node {$(1,1)$};
            \draw (7.5, -0.75) node {$\mathbf{5}$};

       \foreach \i in {6,7,8}{
                 \foreach \j in {0,1,2}{

              \draw (\i,\j) rectangle (\i+1,\j+1);
              }
              }


          \foreach \i in {13}{
                 \foreach \j in {0,2}{

      \fill [pattern= horizontal lines, pattern color =blue]  (\i,\j) rectangle (\i+1,\j+1);
              }
              }

	      \draw (13.5, 3.75) node {$(0,1)$};
            \draw (13.5, -0.75) node {$\mathbf{7}$};

     \foreach \i in {12,13,14}{
                 \foreach \j in {0,1,2}{

              \draw  (\i,\j) rectangle (\i+1,\j+1);
              }
              }

		
  			 \foreach \i in {18,20}{
                 \foreach \j in {1}{

              \fill [pattern= horizontal lines, pattern color =blue]  (\i,\j) rectangle (\i+1,\j+1);
              }
              }

		    \draw (19.5, 3.75) node {$(1,0)$};
              \draw (19.5, -0.75) node {$\mathbf{7}$};

     \foreach \i in {18,19,20}{
                 \foreach \j in {0,1,2}{

              \draw  (\i,\j) rectangle (\i+1,\j+1);
              }
              }

                    \end{tikzpicture}}
\caption{The four possible configurations of $\mathcal{A}_3$. The white squares represent the points of $\mathcal{A}_3(a,b)$; the blue (horizontally lined) squares, points of $\mathcal{B}_3(a,b)$, have both coordinates even. On top, we have noted~$(a,b)$ modulo 2, and below each configuration, we have written 
 the corresponding value of $\Phi_3(a,b)$, ranging from 5 to 8.}
\label{fig:cuadrados3}
\end{figure}

We depict the four possible configurations of $\mathcal{A}_3$ and $\mathcal{B}_3$ in Figure \ref{fig:cuadrados3}.

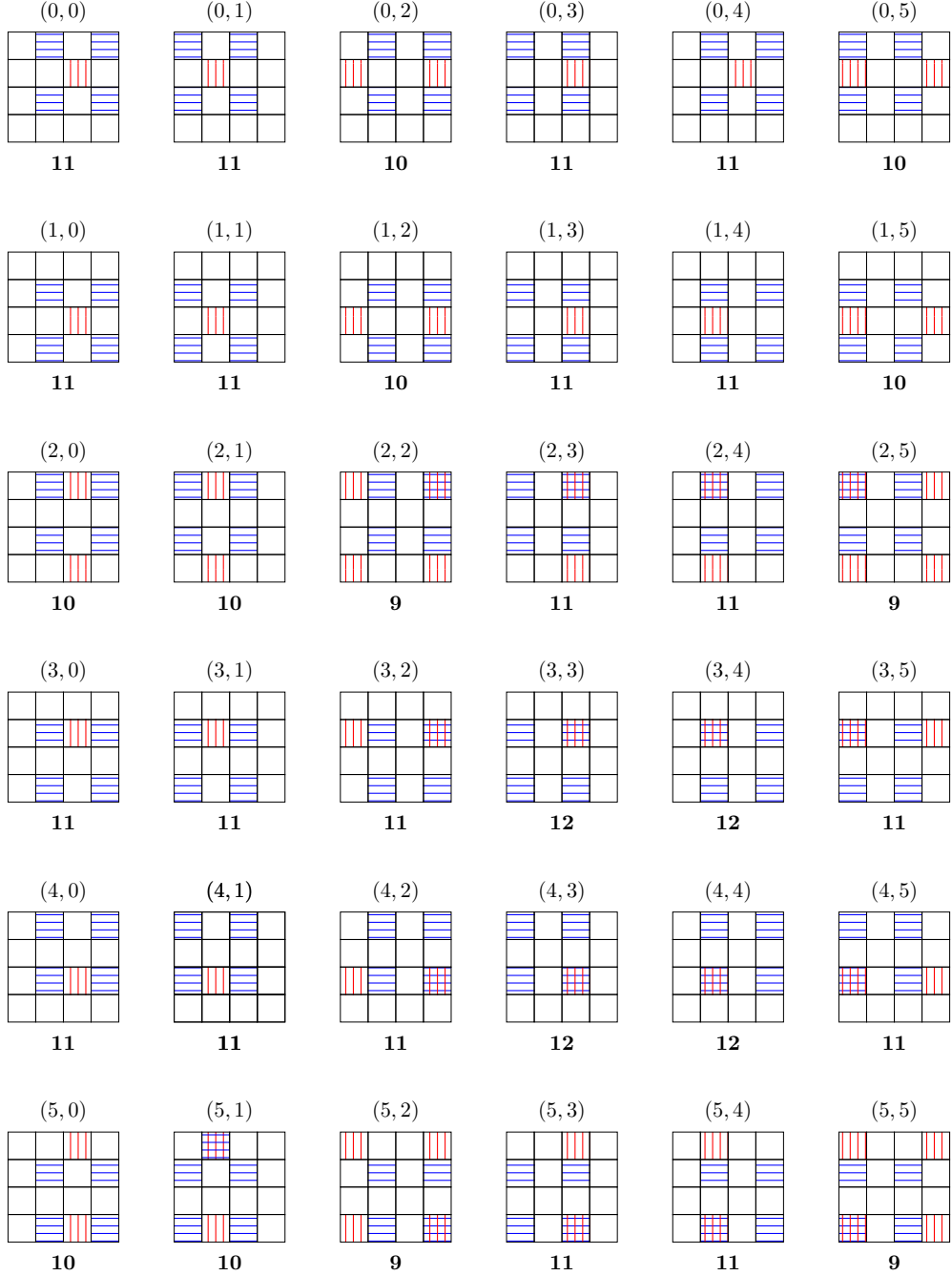
\begin{figure}[h]
\centering\resizebox{13cm}{!}{\begin{tikzpicture}[scale= 0.5]



          \foreach \i in {1,3}{
                 \foreach \j in {0,2}{

              \fill [pattern= horizontal lines, pattern color =blue]  (\i,\j) rectangle (\i+1,\j+1);
              }
              }

     \foreach \i in {2}{
                 \foreach \j in {0,3}{

              \fill [pattern= vertical lines, pattern color =red]  (\i,\j) rectangle (\i+1,\j+1);
              }
              }

            \draw (2, 4.75) node {$(5,0)$};
            \draw (2, -0.75) node {$\mathbf{10}$};

                \foreach \i in {0,1,2,3}{
                 \foreach \j in {0,1,2,3}{

              \draw  (\i,\j) rectangle (\i+1,\j+1);
              }
              }


  			 \foreach \i in {6,8}{
                 \foreach \j in {0,2}{

              \fill [pattern= horizontal lines, pattern color =blue]  (\i,\j) rectangle (\i+1,\j+1);
              }
              }

              \foreach \i in {7}{
                 \foreach \j in {0,3}{

              \fill [pattern= vertical lines, pattern color =red]  (\i,\j) rectangle (\i+1,\j+1);
              }
              }

 		\fill [pattern= horizontal lines, pattern color =blue]  (7,3)  rectangle(8,4);

           \draw (8, 4.75) node {$(5,1)$};
            \draw (8, -0.75) node {$\mathbf{10}$};

       \foreach \i in {6,7,8,9}{
                 \foreach \j in {0,1,2,3}{

              \draw (\i,\j) rectangle (\i+1,\j+1);
              }
              }


          \foreach \i in {13,15}{
                 \foreach \j in {0,2}{

      \fill [pattern= horizontal lines, pattern color =blue]  (\i,\j) rectangle (\i+1,\j+1);
              }
              }

     \foreach \i in {12,15}{
                 \foreach \j in {0,3}{

              \fill [pattern= vertical lines, pattern color =red]  (\i,\j) rectangle (\i+1,\j+1);
              }
              }

	      \draw (14, 4.75) node {$(5,2)$};
            \draw (14, -0.75) node {$\mathbf{9}$};

     \foreach \i in {12,13,14,15}{
                 \foreach \j in {0,1,2,3}{

              \draw  (\i,\j) rectangle (\i+1,\j+1);
              }
              }

		
  			 \foreach \i in {18,20}{
                 \foreach \j in {0,2}{

              \fill [pattern= horizontal lines, pattern color =blue]  (\i,\j) rectangle (\i+1,\j+1);
              }
              }

              \foreach \i in {20}{
                 \foreach \j in {0,3}{

              \fill [pattern= vertical lines, pattern color =red]  (\i,\j) rectangle (\i+1,\j+1);
              }
              }


		    \draw (20, 4.75) node {$(5,3)$};
              \draw (20, -0.75) node {$\mathbf{11}$};

     \foreach \i in {18,19,20,21}{
                 \foreach \j in {0,1,2,3}{

              \draw  (\i,\j) rectangle (\i+1,\j+1);
              }
              }


          \foreach \i in {25,27}{
                 \foreach \j in {0,2}{

		  \fill [pattern= horizontal lines, pattern color =blue]  (\i,\j) rectangle (\i+1,\j+1);
              }
              }

     \foreach \i in {25}{
                 \foreach \j in {0,3}{

              \fill [pattern= vertical lines, pattern color =red]  (\i,\j) rectangle (\i+1,\j+1);
              }
              }

		  \draw (26, 4.75) node {$(5,4)$};
              \draw (26, -0.75) node {$\mathbf{11}$};

   	 \foreach \i in {24,25,26,27}{
                 \foreach \j in {0,1,2,3}{

              \draw  (\i,\j) rectangle (\i+1,\j+1);
              }
              }

		
  			 \foreach \i in {30,32}{
                 \foreach \j in {0,2}{

              \fill [pattern= horizontal lines, pattern color =blue]  (\i,\j) rectangle (\i+1,\j+1);
              }
              }

              \foreach \i in {30,33}{
                 \foreach \j in {0,3}{

              \fill [pattern= vertical lines, pattern color =red]  (\i,\j) rectangle (\i+1,\j+1);
              }
              }


		  \draw (32, 4.75) node {$(5,5)$};
              \draw (32, -0.75) node {$\mathbf{9}$};

	   \foreach \i in {30,31,32,33}{
                 \foreach \j in {0,1,2,3}{

              \draw  (\i,\j) rectangle (\i+1,\j+1);
              }
              }


          \foreach \i in {1,3}{
                 \foreach \j in {9,11}{

              \fill [pattern= horizontal lines, pattern color =blue]  (\i,\j) rectangle (\i+1,\j+1);
              }
              }

              \fill [pattern= vertical lines, pattern color =red]  (2,9) rectangle (3,10);

            \draw (2, 12.75) node {$(4,0)$};
            \draw (2, 7.25) node {$\mathbf{11}$};

                \foreach \i in {0,1,2,3}{
                 \foreach \j in {8,9,10,11}{

              \draw  (\i,\j) rectangle (\i+1,\j+1);
              }
              }


 \foreach \i in {6,8}{
                 \foreach \j in {9,11}{

              \fill [pattern= horizontal lines, pattern color =blue]  (\i,\j) rectangle (\i+1,\j+1);
              }
              }
 		
             \fill [pattern= vertical lines, pattern color =red]  (7,9) rectangle (8,10);

           \draw (8, 12.75) node {$(4,1)$};
            \draw (8, 7.25) node {$\mathbf{11}$};

       \foreach \i in {6,7,8,9}{
                 \foreach \j in {8,9,10,11}{

              \draw (\i,\j) rectangle (\i+1,\j+1);
              }
              }


		 \foreach \i in {13,15}{
                 \foreach \j in {9,11}{

              \fill [pattern= horizontal lines, pattern color =blue]  (\i,\j) rectangle (\i+1,\j+1);
              }
              }

 		 \foreach \i in {12,15}{
                 \foreach \j in {9}{
              \fill [pattern= vertical lines, pattern color =red] (\i,\j) rectangle (\i+1,\j+1);

		     }
              }

	      \draw (14, 12.75) node {$(4,2)$};
            \draw (14, 7.25) node {$\mathbf{11}$};

     \foreach \i in {12,13,14,15}{
                 \foreach \j in {8,9,10,11}{

              \draw  (\i,\j) rectangle (\i+1,\j+1);
              }
              }


		\foreach \i in {18,20}{
                 \foreach \j in {9,11}{

              \fill [pattern= horizontal lines, pattern color =blue]  (\i,\j) rectangle (\i+1,\j+1);
              }
              }
 		
             \fill [pattern= vertical lines, pattern color =red]  (20,9) rectangle (21,10);

           \draw (8, 12.75) node {$(4,1)$};
            \draw (8, 7.25) node {$\mathbf{11}$};

       \foreach \i in {6,7,8,9}{
                 \foreach \j in {8,9,10,11}{

              \draw (\i,\j) rectangle (\i+1,\j+1);
              }
              }

		    \draw (20, 12.75) node {$(4,3)$};
              \draw (20, 7.25) node {$\mathbf{12}$};

     \foreach \i in {18,19,20,21}{
                 \foreach \j in {8,9,10,11}{

              \draw  (\i,\j) rectangle (\i+1,\j+1);
              }
              }

		\foreach \i in {25,27}{
                 \foreach \j in {9,11}{

              \fill [pattern= horizontal lines, pattern color =blue]  (\i,\j) rectangle (\i+1,\j+1);
              }
              }

               \fill [pattern= vertical lines, pattern color =red]  (25,9) rectangle (26,10);

		  \draw (26, 12.75) node {$(4,4)$};
              \draw (26, 7.25) node {$\mathbf{12}$};

   	 \foreach \i in {24,25,26,27}{
                 \foreach \j in {8,9,10,11}{

              \draw  (\i,\j) rectangle (\i+1,\j+1);
              }
              }

		\foreach \i in {30,32}{
                 \foreach \j in {9,11}{

              \fill [pattern= horizontal lines, pattern color =blue]  (\i,\j) rectangle (\i+1,\j+1);
              }
              }

              		 \foreach \i in {30,33}{
                 \foreach \j in {9}{
              \fill [pattern= vertical lines, pattern color =red] (\i,\j) rectangle (\i+1,\j+1);

		     }
              }

		  \draw (32, 12.75) node {$(4,5)$};
              \draw (32, 7.25) node {$\mathbf{11}$};

	   \foreach \i in {30,31,32,33}{
                 \foreach \j in {8,9,10,11}{

              \draw  (\i,\j) rectangle (\i+1,\j+1);
              }
              }


 \foreach \i in {1,3}{
                 \foreach \j in {16,18}{

              \fill [pattern= horizontal lines, pattern color =blue]  (\i,\j) rectangle (\i+1,\j+1);
              }
              }

              \fill [pattern= vertical lines, pattern color =red]  (2, 18) rectangle (3,19);

            \draw (2, 20.75) node {$(3,0)$};
            \draw (2, 15.25) node {$\mathbf{11}$};

                \foreach \i in {0,1,2,3}{
                 \foreach \j in {16,17,18,19}{

              \draw  (\i,\j) rectangle (\i+1,\j+1);
              }
              }


 \foreach \i in {6,8}{
                 \foreach \j in {16,18}{

              \fill [pattern= horizontal lines, pattern color =blue]  (\i,\j) rectangle (\i+1,\j+1);
              }
              }
 		
             \fill [pattern= vertical lines, pattern color =red]  (7,18) rectangle (8,19);

           \draw (8, 20.75) node {$(3,1)$};
            \draw (8, 15.25) node {$\mathbf{11}$};

       \foreach \i in {6,7,8,9}{
                 \foreach \j in {16,17,18,19}{

              \draw (\i,\j) rectangle (\i+1,\j+1);
              }
              }


		 \foreach \i in {13,15}{
                 \foreach \j in {16,18}{

              \fill [pattern= horizontal lines, pattern color =blue]  (\i,\j) rectangle (\i+1,\j+1);
              }
              }
	
  \foreach \i in {12,15}{
                 \foreach \j in {18}{

              \fill [pattern= vertical lines, pattern color =red]  (\i,\j) rectangle (\i+1,\j+1);
              }
              }
	      \draw (14, 20.75) node {$(3,2)$};
            \draw (14, 15.25) node {$\mathbf{11}$};

     \foreach \i in {12,13,14,15}{
                 \foreach \j in  {16,17,18,19}{

              \draw  (\i,\j) rectangle (\i+1,\j+1);
              }
              }


		 \foreach \i in {18,20}{
                 \foreach \j in {16,18}{

              \fill [pattern= horizontal lines, pattern color =blue]  (\i,\j) rectangle (\i+1,\j+1);
              }
              }

		 \fill [pattern= vertical lines, pattern color =red]  (20,18) rectangle (21,19);

		    \draw (20, 20.75) node {$(3,3)$};
              \draw (20, 15.25) node {$\mathbf{12}$};

     \foreach \i in {18,19,20,21}{
                 \foreach \j in {16,17,18,19}{

              \draw  (\i,\j) rectangle (\i+1,\j+1);
              }
              }

		 \foreach \i in {25,27}{
                 \foreach \j in {16,18}{

              \fill [pattern= horizontal lines, pattern color =blue]  (\i,\j) rectangle (\i+1,\j+1);
              }
              }

            \fill [pattern= vertical lines, pattern color =red]  (25,18) rectangle (26,19);

		  \draw (26, 20.75) node {$(3,4)$};
              \draw (26, 15.25) node {$\mathbf{12}$};

   	 \foreach \i in {24,25,26,27}{
                 \foreach \j in {16,17,18,19}{

              \draw  (\i,\j) rectangle (\i+1,\j+1);
              }
              }

		 \foreach \i in {30,32}{
                 \foreach \j in {16,18}{

              \fill [pattern= horizontal lines, pattern color =blue]  (\i,\j) rectangle (\i+1,\j+1);
              }
              }

		 \foreach \i in {30,33}{
                 \foreach \j in {18}{

              \fill [pattern= vertical lines, pattern color =red]  (\i,\j) rectangle (\i+1,\j+1);
              }
              }
		  \draw (32, 20.75) node {$(3,5)$};
              \draw (32, 15.25) node {$\mathbf{11}$};

	   \foreach \i in {30,31,32,33}{
                 \foreach \j in {16,17,18,19}{

              \draw  (\i,\j) rectangle (\i+1,\j+1);
              }
              }


      \foreach \i in {1,3}{
                 \foreach \j in {25,27}{

              \fill [pattern= horizontal lines, pattern color =blue]  (\i,\j) rectangle (\i+1,\j+1);
              }
              }
            \foreach \i in {2}{
                 \foreach \j in {24,27}{

              \fill [pattern= vertical lines, pattern color =red]  (\i,\j) rectangle (\i+1,\j+1);
              }
              }


            \draw (2, 28.75) node {$(2,0)$};
            \draw (2, 23.25) node {$\mathbf{10}$};

                \foreach \i in {0,1,2,3}{
                 \foreach \j in {24,25,26,27}{

              \draw  (\i,\j) rectangle (\i+1,\j+1);
              }
              }


		 \foreach \i in {6,8}{
                 \foreach \j in {25,27}{

              \fill [pattern= horizontal lines, pattern color =blue]  (\i,\j) rectangle (\i+1,\j+1);
              }
              }
 		
              \foreach \i in {7}{
                 \foreach \j in {24,27}{

              \fill [pattern= vertical lines, pattern color =red]  (\i,\j) rectangle (\i+1,\j+1);
              }
              }

           \draw (8, 28.75) node {$(2,1)$};
            \draw (8, 23.25) node {$\mathbf{10}$};

       \foreach \i in {6,7,8,9}{
                 \foreach \j in  {24,25,26,27}{

              \draw (\i,\j) rectangle (\i+1,\j+1);
              }
              }


		 \foreach \i in {13,15}{
                 \foreach \j in {25,27}{

              \fill [pattern= horizontal lines, pattern color =blue]  (\i,\j) rectangle (\i+1,\j+1);
              }
              }
            \foreach \i in {12,15}{
                 \foreach \j in {24,27}{

              \fill [pattern= vertical lines, pattern color =red]  (\i,\j) rectangle (\i+1,\j+1);
              }
              }

	      \draw (14, 28.75) node {$(2,2)$};
            \draw (14, 23.25) node {$\mathbf{9}$};

     \foreach \i in {12,13,14,15}{
                 \foreach \j in  {24,25,26,27}{

              \draw  (\i,\j) rectangle (\i+1,\j+1);
              }
              }


  		 \foreach \i in {18,20}{
                 \foreach \j in {25,27}{

              \fill [pattern= horizontal lines, pattern color =blue]  (\i,\j) rectangle (\i+1,\j+1);
              }
              }
 		
              \foreach \i in {20}{
                 \foreach \j in {24,27}{

              \fill [pattern= vertical lines, pattern color =red]  (\i,\j) rectangle (\i+1,\j+1);
              }
              }

		    \draw (20, 28.75) node {$(2,3)$};
              \draw (20, 23.25) node {$\mathbf{11}$};

     \foreach \i in {18,19,20,21}{
                 \foreach \j in  {24,25,26,27}{

              \draw  (\i,\j) rectangle (\i+1,\j+1);
              }
              }


				 \foreach \i in {25,27}{
                 \foreach \j in {25,27}{

              \fill [pattern= horizontal lines, pattern color =blue]  (\i,\j) rectangle (\i+1,\j+1);
              }
              }

		   \foreach \i in {25}{
                 \foreach \j in {24,27}{

              \fill [pattern= vertical lines, pattern color =red]  (\i,\j) rectangle (\i+1,\j+1);
              }
              }

		  \draw (26, 28.75) node {$(2,4)$};
              \draw (26, 23.25) node {$\mathbf{11}$};

   	 \foreach \i in {24,25,26,27}{
                 \foreach \j in  {24,25,26,27}{

              \draw  (\i,\j) rectangle (\i+1,\j+1);
              }
              }


  		 \foreach \i in {30,32}{
                 \foreach \j in {25,27}{

              \fill [pattern= horizontal lines, pattern color =blue]  (\i,\j) rectangle (\i+1,\j+1);
              }
              }
 		
              \foreach \i in {30,33}{
                 \foreach \j in {24,27}{

              \fill [pattern= vertical lines, pattern color =red]  (\i,\j) rectangle (\i+1,\j+1);
              }
              }

		  \draw (32, 28.75) node {$(2,5)$};
              \draw (32, 23.25) node {$\mathbf{9}$};

	   \foreach \i in {30,31,32,33}{
                 \foreach \j in  {24,25,26,27}{

              \draw  (\i,\j) rectangle (\i+1,\j+1);
              }
              }


               		 \foreach \i in {1,3}{
                 \foreach \j in {32,34}{

              \fill [pattern= horizontal lines, pattern color =blue]  (\i,\j) rectangle (\i+1,\j+1);
              }
              }

              \fill [pattern= vertical lines, pattern color =red]  (2,33) rectangle (3,34);

            \draw (2, 36.75) node {$(1,0)$};
            \draw (2, 31.25) node {$\mathbf{11}$};

                \foreach \i in {0,1,2,3}{
                 \foreach \j in {32,33,34,35}{

              \draw  (\i,\j) rectangle (\i+1,\j+1);
              }
              }


 		  \foreach \i in {6,8}{
                 \foreach \j in {32,34}{

              \fill [pattern= horizontal lines, pattern color =blue]  (\i,\j) rectangle (\i+1,\j+1);
              }
              }

 		       \fill [pattern= vertical lines, pattern color =red]  (7,33) rectangle (8,34);

           \draw (8, 36.75) node {$(1,1)$};
            \draw (8, 31.25) node {$\mathbf{11}$};

       \foreach \i in {6,7,8,9}{
                 \foreach \j in  {32,33,34,35}{

              \draw (\i,\j) rectangle (\i+1,\j+1);
              }
              }


 		 \foreach \i in {13,15}{
                 \foreach \j in {32,34}{

              \fill [pattern= horizontal lines, pattern color =blue]  (\i,\j) rectangle (\i+1,\j+1);
              }
              }

         \foreach \i in {12,15}{
                 \foreach \j in {33}{

              \fill [pattern= vertical lines, pattern color =red]  (\i,\j) rectangle (\i+1,\j+1);
              }
              }

	      \draw (14, 36.75) node {$(1,2)$};
            \draw (14, 31.25) node {$\mathbf{10}$};

     \foreach \i in {12,13,14,15}{
                 \foreach \j in  {32,33,34,35}{

              \draw  (\i,\j) rectangle (\i+1,\j+1);
              }
              }


		\foreach \i in {18,20}{
                 \foreach \j in {32,34}{

              \fill [pattern= horizontal lines, pattern color =blue]  (\i,\j) rectangle (\i+1,\j+1);
              }
              }

 		       \fill [pattern= vertical lines, pattern color =red]  (20,33) rectangle (21,34);

		    \draw (20, 36.75) node {$(1,3)$};
              \draw (20, 31.25) node {$\mathbf{11}$};

     \foreach \i in {18,19,20,21}{
                 \foreach \j in {32,33,34,35}{

              \draw  (\i,\j) rectangle (\i+1,\j+1);
              }
              }


 		 \foreach \i in {25,27}{
                 \foreach \j in {32,34}{

              \fill [pattern= horizontal lines, pattern color =blue]  (\i,\j) rectangle (\i+1,\j+1);
              }
              }

 		       \fill [pattern= vertical lines, pattern color =red]  (25,33) rectangle (26,34);

		  \draw (26, 36.75) node {$(1,4)$};
              \draw (26, 31.25) node {$\mathbf{11}$};

   	 \foreach \i in {24,25,26,27}{
                 \foreach \j in {32,33,34,35}{

              \draw  (\i,\j) rectangle (\i+1,\j+1);
              }
              }


		\foreach \i in {30,32}{
                 \foreach \j in {32,34}{

              \fill [pattern= horizontal lines, pattern color =blue]  (\i,\j) rectangle (\i+1,\j+1);
              }
              }
		  \foreach \i in {30,33}{
                 \foreach \j in {33}{

              \fill [pattern= vertical lines, pattern color =red]  (\i,\j) rectangle (\i+1,\j+1);
              }
              }

		  \draw (32, 36.75) node {$(1,5)$};
              \draw (32, 31.25) node {$\mathbf{10}$};

	   \foreach \i in {30,31,32,33}{
                 \foreach \j in {32,33,34,35}{

              \draw  (\i,\j) rectangle (\i+1,\j+1);
              }
              }



   \foreach \i in {1,3}{
                 \foreach \j in {41,43}{

              \fill [pattern= horizontal lines, pattern color =blue]  (\i,\j) rectangle (\i+1,\j+1);
              }
              }
 		
             \fill [pattern= vertical lines, pattern color =red]  (2,42) rectangle (3,43);

            \draw (2, 44.75) node {$(0,0)$};
            \draw (2, 39.25) node {$\mathbf{11}$};

                \foreach \i in {0,1,2,3}{
                 \foreach \j in {40,41,42,43}{

              \draw  (\i,\j) rectangle (\i+1,\j+1);
              }
              }


     \foreach \i in {6,8}{
                 \foreach \j in {41,43}{

              \fill [pattern= horizontal lines, pattern color =blue]  (\i,\j) rectangle (\i+1,\j+1);
              }
              }
 		
             \fill [pattern= vertical lines, pattern color =red]  (7,42) rectangle (8,43);

           \draw (8, 44.75) node {$(0,1)$};
            \draw (8, 39.25) node {$\mathbf{11}$};

       \foreach \i in {6,7,8,9}{
                 \foreach \j in {40,41,42,43}{

              \draw (\i,\j) rectangle (\i+1,\j+1);
              }
              }


  		 \foreach \i in {13,15}{
                 \foreach \j in {41,43}{

              \fill [pattern= horizontal lines, pattern color =blue]  (\i,\j) rectangle (\i+1,\j+1);
              }
              }

 				  \foreach \i in {12,15}{
                 \foreach \j in {42}{

              \fill [pattern= vertical lines, pattern color =red]  (\i,\j) rectangle (\i+1,\j+1);
              }
              }

	      \draw (14, 44.75) node {$(0,2)$};
            \draw (14, 39.25) node {$\mathbf{10}$};

     \foreach \i in {12,13,14,15}{
                 \foreach \j in {40,41,42,43}{

              \draw  (\i,\j) rectangle (\i+1,\j+1);
              }
              }


     \foreach \i in {18,20}{
                 \foreach \j in {41,43}{

              \fill [pattern= horizontal lines, pattern color =blue]  (\i,\j) rectangle (\i+1,\j+1);
              }
              }
 		
             \fill [pattern= vertical lines, pattern color =red]  (20,42) rectangle (21,43);

		    \draw (20, 44.75) node {$(0,3)$};
              \draw (20, 39.25) node {$\mathbf{11}$};

     \foreach \i in {18,19,20,21}{
                 \foreach \j in {40,41,42,43}{

              \draw  (\i,\j) rectangle (\i+1,\j+1);
              }
              }


  			 \foreach \i in {25,27}{
                 \foreach \j in {41,43}{

              \fill [pattern= horizontal lines, pattern color =blue]  (\i,\j) rectangle (\i+1,\j+1);
              }
              }
 		
             \fill [pattern= vertical lines, pattern color =red]  (26,42) rectangle (27,43);

		  \draw (26, 44.75) node {$(0,4)$};
              \draw (26, 39.25) node {$\mathbf{11}$};

   	 \foreach \i in {24,25,26,27}{
                 \foreach \j in {40,41,42,43}{

              \draw  (\i,\j) rectangle (\i+1,\j+1);
              }
              }


		 \foreach \i in {30,32}{
                 \foreach \j in {41,43}{

              \fill [pattern= horizontal lines, pattern color =blue]  (\i,\j) rectangle (\i+1,\j+1);
              }
              }

 		    \foreach \i in {30, 33}{
                 \foreach \j in {42}{

              \fill [pattern= vertical lines, pattern color =red]  (\i,\j) rectangle (\i+1,\j+1);
              }
              }

		  \draw (32, 44.75) node {$(0,5)$};
              \draw (32, 39.25) node {$\mathbf{10}$};

	   \foreach \i in {30,31,32,33}{
                 \foreach \j in {40,41,42,43}{

              \draw  (\i,\j) rectangle (\i+1,\j+1);
              }
              }

                    \end{tikzpicture}}
          \caption{The 36 possible configurations of $\mathcal{A}_4$, labelled  with the values of $(a,b)$ modulo~6, and the corresponding values of $\Phi_4(a,b)$, that in this case range from~9 to~12.
          As before, the white squares represent the points of $\mathcal{A}_4(a,b)$, and the blue (horizontally lined) squares correspond to points with both coordinates even; but now the red squares (vertically lined) have both coordinates divisible by~3. Some squares, of course, belong to both categories.}
         \label{fig:cuadrados4}
         \end{figure}

For the case $M=4$, with $P_4=6$, Figure \ref{fig:cuadrados4} displays the 36 possible configurations.

For fixed $s$, such that $0 \le s\le M^2$, we denote with $\xi
 (M,s)$ the (arithmetic) average  of the binomial coeficientes $\Phi_M(u,v)\choose s$:
\begin{equation}\label{eq:formula of gammas}
\xi(M,s)=\frac{1}{P_M^2} \sum_{0\le u,v<P_M} \binom{\Phi_M(u,v)}{s} , \quad\mbox{for $0\le s \le M^2$}.
\end{equation}
Of course, $\xi(M,0)=1$. As we shall see later, in \eqref{eq:formula gamma(M,1)}, the value of
$\xi(M,1)$, which is the average value of $\Phi_M(u,v)$, is 
\[
\xi(M,1)=\frac{1}{P_M^2} \sum_{0\le u,v<P_M} \Phi_M(u,v)= M^2 \prod_{p<M} \Big(1-\frac{1}{p^2}\Big).
\]

\subsection{Distribution of the  window coprime counting variable}\label{sec:distr of Z and Z*}

Fix $M\ge 1$. We use now the splitting of Section \ref{sec:splitting} to simplify the definition of $Z_M(a,b)$ from the expression~\eqref{eq:def of ZM}  
to
\begin{align} \nonumber
Z_M(a,b)&=\sum_{(k,l) \in \K_M} V(a + k,b + l)
=\sum_{(k,l) \in \B_M(a,b)} V(a + k,b+l)+ \sum_{(k,l) \in \A_M(a,b)} V(a + k,b + l)
\\ \nonumber
&
=\sum_{(k,l) \in \A_M(a,b)} V(a + k,b + l)
=\sum_{(k,l) \in  \A_M(a,b)} \ \prod_{p} (1- I_p(a  +   k)\,I_p(b +   l))
\\ \label{eq:def de ZM como suma de indicadores en A}
&=\sum_{(k,l) \in  \A_M(a,b)} \ \prod_{p\ge M} (1- I_p(a +  k)\,I_p(b +  l)).
\end{align}
Here, we have used \eqref{eq:formula V(a,b)} and the very definitions of $\mathcal{A}_M(a,b)$ and $\mathcal{B}_M(a,b)$.
Notice also how this shows that the maximum possible value of $Z_M(a,b)$ is $\Phi(a,b)$, and not $M^2$.

\subsubsection{Conditioning upon residues modulo $P_M$}

We now fix a pair $(u, v)$ of  residues modulo~$ P_M$, with $0 \le u,v  < 
 P_M$,

Fix $n\ge 1$. Recall, from \eqref{eq:def of RM}, that $R_M(a,b)$ denotes the pair $(r_{P_M}(a),r_{P_M}(b))$ of  residues of $(a,b)$ modulo $P_M$. We are going to condition upon
\begin{equation}\label{eq:def Omega_n}
\Omega_n(u, v)=\{(a,b)\in \N_n^2: R_M(a,b)=(u, v)\}.
\end{equation}
There are a total of $P_M^2$ different $\Omega_n(u,v)$,  which form a partition of $\N_n^2$.

Using Lemma \ref{lemma:asympotic independence}, we see  that
\begin{equation}\label{eq:lim prob Omega_n} \lim_{n\to \infty} \P_n(\Omega_n(u,v))=\frac{1}{P_M^2} , \quad \mbox{for each $0 \le u,v <P_M$}\cdot
\end{equation}

For all $(a,b) \in \Omega_n(u,v)$ we have that $\mathcal{A}_M(a,b)=\A_M(u,v)$, since $R_M(a,b)=(u,v)$,  that is, $a\equiv u \mod p$ and $b\equiv v \mod p$ for all prime $p>M$. Consequently, we have that  $\Phi_M(a,b)=\Phi_M(u,v)$ if $(a,b) \in \Omega_n(u,v)$.

\medskip

We apply now Lemma \ref{lemma:schuette-nesbitt}. Take, in the notation used there, $\Omega_n(u,v)$ as $\Omega$,  the function~$Z_M$ as the counting function $C$ and $\Phi_M(u,v)$ as $t$. Then, for any $(a,b)\in \Omega_n(u,v)$ and for $r$ such that $0 \le r \le \Phi_M(u,v)$, we have that
\begin{equation}\label{eq:indicadora en Omegan1}
\begin{aligned}
\uno\nolimits_{\{Z_M=r\}}(a,b)
&=\sum_{s=r}^{\Phi_M(u,v)} (-1)^{s-r} \binom{s}{r} \sum_{\substack{H \subset \A_M(u,v) \\
|H|=s}}  \prod_{(k,l)\in H}  \prod_{p\ge M} (1-I_p(a +  k)\,I_p(b + l)).
\\
&=\sum_{s=r}^{\Phi_M(u,v)} (-1)^{s-r} \binom{s}{r} \sum_{\substack{H \subset \A_M(u,v) \\
|H|=s}}  \prod_{p\ge M} \prod_{(k,l)\in H}   (1-I_p(a +  k)\,I_p(b + l)).
\end{aligned}\end{equation}

Now observe that, if $(k,l)\neq (k', l') \in \K_M$, then
\[
\begin{aligned}
\big[1-I_p(a +  k) &\,I_p(b+l)\big]\,\big[1-I_p(a +  k')\,I_p(b +  l')\big]
\\
&\qquad\qquad =1-I_p(a +  k)\,I_p(b +  l)-I_p(a +  k')\,I_p(b +  l') ,\end{aligned}
\]
because the term 
 $I_p(a +  k)\,I_p(b +  l)\,I_p(a +  k')\,I_p(b +  l')$ 
  vanishes.
This is so since if $k \neq k'$, then $I_p(a +  k)\,I_p(a +  k')=0$, because the prime $p\ge M$ does not divide $k - k'$, and analogously, if $l\neq l'$, then $I_p(b +  l)\,I_p(b +  l')=0$.

Using this observation, we may rewrite \eqref{eq:indicadora en Omegan1} as follows: for any $(a,b)\in \Omega_n(u,v)$ and for~$r$ such that $0 \le r \le \Phi(u,v)$,
\begin{equation}\label{eq:indicadora en Omegan2}
\uno\nolimits_{\{Z_M=r\}}(a,b)
=\sum_{s=r}^{\Phi_M(u,v)} (-1)^{s-r}\! \binom{s}{r} \sum_{\substack{H \subset \A_M(u,v) \\
|H|=s}} \,   \prod_{p\ge M}\!\! \Big(1-\sum_{(k,l)\in H}\!\! I_p(a + k)\,I_p(b +  l)\Big) .
\end{equation}

Now, for each $p \ge M$, the function $(a,b)\mapsto 1-\sum_{(k,l)\in H} I_p(a +  k)\,I_p(b +  l)$ takes only the values 0 and 1. For if $p \ge M$, then $p$ may divide at most one $a+k$ with $1\le k \le M$ (and also at most one $b+l$ with $1 \le l \le M$). Thus this function  is the indicator function of a certain subset in  $\N^2$, which we denote by  $B_H^{(p)}$:
\begin{equation}\label{eq:def de Bhp}
\uno\nolimits_{B_H^{(p)}}(a,b)=1-\sum_{(k,l)\in H} I_p(a +  k)\,I_p(b +  l).
\end{equation}
With this new notation, we may finally rewrite, for any $(a,b)\in \Omega_n(u,v)$ and for $r$ such that $0 \le r \le \Phi(u,v)$,
\begin{align} \nonumber
\uno\nolimits_{\{Z_M=r\}}(a,b)
&
=\sum_{s=r}^{\Phi_M(u,v)} (-1)^{s-r} \binom{s}{r} \sum_{\substack{H \subset \A_M(u,v) \\
|H|=s}} \,   \prod_{p\ge M} \uno\nolimits_{B_H^{(p)}}(a,b)
\\ \label{eq:indicadora en Omegan3}
&=\sum_{s=r}^{\Phi_M(u,v)} (-1)^{s-r} \binom{s}{r} \sum_{\substack{H \subset \A_M(u,v) \\
|H|=s}} \,   \uno\nolimits_{\{\bigcap_{p\ge M} B_H^{(p)}\}}(a,b).
\end{align}

Observe that for each $(a,b) \in \N_{n}^2$, the above products (or the above intersections) are actually products/intersections of finitely  many terms/sets, since for $p>n+M$, we have $I_p(a +  k)\,I_p(b +  l)=0$, for any $(k,l)\in \K_M$.

Regarding these $B_H^{(p)}$, we have the following key lemma.

\begin{lemma}\label{lemma:a la caibach} For any $H \subset \mathcal{A}_M(u,v)$ such that $|H|=s$, we have that
\[
\lim_{n \to \infty} \P_n\Big(\bigcap_{p \ge M}B_H^{(p)}\,|\,\Omega_n(u,v)\Big)=\prod_{p\ge M} \Big(1-\frac{s}{p^2}\Big).
\]
\end{lemma}
\begin{proof} This argument is modeled upon the proof of Theorem~5 in~\cite{CaiBach}.
Fix $M\ge 1$ and a pair~$(u,v)$ such that $0\le u,v<P_M$. Consider the sets $\Omega_n(u,v)$ and $\mathcal{A}_M(u,v)$ defined in~\eqref{eq:def Omega_n} and~\eqref{eq:def A_m}, respectively.

Now fix a subset~$H$ of $\mathcal{A}_M(u,v)$ of size $|H|=s$. Observe first that for $p \ge M$, and using the definition of $B_H^{(p)}$ in~\eqref{eq:def de Bhp},
\begin{align}\nonumber
\P_n\big(\N_n^2\setminus B_H^{(p)} &\,|\,\Omega_n(u,v)\big)
=\E_n \Big(\sum_{(k,l)\in H} I_p(\cdot +k)\,I_p(\cdot+l) \,|\, \Omega_n(u,v)\Big)
\\
\nonumber
&=\sum_{(k,l)\in H} \E_n(I_p(\cdot +k)\,I_p(\cdot+l) \,|\, \Omega_n(u,v))
\\
&=\sum_{(k,l)\in H}\P_n\big(\{a,b\in\mathbb{N}^2_n: a+k\equiv 0 \textup{ mod }p, b+l\equiv 0 \textup{ mod }p\}\,|\, \Omega_n(u,v)\big).
 \label{eq:Pn de Nn menos BHp}
\end{align}
Using that $\Omega_n(u,v)$ is defined in terms of residues modulo the primes $q<M$, while the prime $p$ defining $B_H^{(p)}$ is $p\ge M$,  Lemma \ref{lemma:asympotic independence} gives that
\[
\lim_{n \to \infty}\P_n\big(\N^2_n\setminus B_H^{(p)}\,|\,\Omega_n(u,v)\big)=\frac{s}{p^2}\cdot
\]
Analogously, for any distinct primes $p_1, \ldots, p_R$, all $\ge M$, we have that
\begin{equation}\label{eq:cota interseccion en primos}
\lim_{n \to \infty}\P_n\Big(\bigcap_{i=1}^R \big(\N_n^2\setminus B_H^{(p_i)}\big)\,|\,\Omega_n(u,v)\Big)=s^R\prod_{i=1}^R \frac{1}{p_i^2}\cdot
\end{equation}

Take $N>M$, and observe that
\[
\P_n\Big(\bigcap_{M\le p \le N}B_H^{(p)}\,|\,\Omega_n(u,v)\Big)
=\P_n\Big(\mathbb{N}^2_n\setminus  \bigcup_{M\le p \le N} (\mathbb{N}^2_n\setminus B_H^{(p)}) \,|\,\Omega_n(u,v)\Big).
\]
From \eqref{eq:cota interseccion en primos} and the inclusion/exclusion principle stated in \eqref{eq:IEprob}, 
 and as already discussed in Lemmas \ref{lemma:cesaro identity} and \ref{lemma:formula ZM}, we deduce that 
\begin{equation}\label{eq:interseccion BHp condic}
\lim_{n \to \infty} \P_n\Big(\bigcap_{M\le p \le N}B_H^{(p)}\,|\,\Omega_n(u,v)\Big)=\prod_{M\le p \le N} \Big(1-\frac{s}{p^2}\Big).
\end{equation}
As this holds for for any $N>M$, we deduce that 
\[
\limsup_{n \to \infty} \P_n\Big(\bigcap_{ p \ge M}B_H^{(p)}\,|\,\Omega_n(u,v)\Big)\le\prod_{p\ge M} \Big(1-\frac{s}{p^2}\Big).
\]

For an inequality with $\liminf$ in the opposite  direction, we argue as follows.
For $k,l\le M$, we have that
\[
\#\{(a,b)\in\mathbb{N}^2_n: a+k\equiv 0 \textup{ mod }p, b+l\equiv 0 \textup{ mod }p\}
\le
\#\{(c,d)\in\mathbb{N}^2_{n+M}: p\mid c \text{ and } p\mid d\},
\]
and so,
\[
\P_n\big(\{(a,b)\in\mathbb{N}^2_n: a+k\equiv 0 \textup{ mod }p, b+l\equiv 0 \textup{ mod }p\}\,|\, \Omega_n(u,v)\big)
\le \frac{1}{n^2}\, \Big\lfloor \frac{n +  M}{p}\Big\rfloor^2 \frac{1}{1/P_M^2} \cdot
\]
 This 
 (rather crude) estimate is enough for our purposes. Now, going back to \eqref{eq:Pn de Nn menos BHp}, we find that
for some constant $C_M$, depending on $M$, and for $p \ge M$,
\begin{equation}\label{eq:cota Pn de Nn menos BHp}
\P_n\big(\N_n^2\setminus B_H^{(p)}\,|\,\Omega_n(u,v)\big)\le C_M \,\frac{s}{p^2}\cdot
\end{equation}
For $N >M$, we have that
\[\begin{aligned}
&\P_n\Big(\bigcap_{M\le p \le N}B_H^{(p)}\,|\,\Omega_n(u,v)\Big)-\P_n\Big(\bigcap_{ p \ge M}B_H^{(p)}\,|\,\Omega_n(u,v)\Big)
\\
&=\P_n\Big(\bigcup_{p \ge M}(\mathbb{N}^2_n\setminus B_H^{(p)})\,|\,\Omega_n(u,v)\Big)
-\P_n\Big(\bigcup_{M\le p \le N}(\mathbb{N}^2_n\setminus B_H^{(p)})\,|\,\Omega_n(u,v)\Big)
\\
&=\P_n\Big(\bigcup_{p >N}\big(\mathbb{N}^2_n\setminus B_H^{(p)}\big)\,|\,\Omega_n(u,v)\Big)
\le \sum_{p >N} \P_n\big((\N_n^2\setminus B_H^{(p)})\,|\,\Omega_n(u,v)\big)
\le C_M \, s\, \sum_{p>N}\frac{1}{p^2} ,\end{aligned}
\]
where \eqref{eq:cota Pn de Nn menos BHp} was used in the last inequality. Thus, we can 
 deduce, using \eqref{eq:interseccion BHp condic}, that, for $N >M$,
\[
\liminf_{n \to \infty} \P_n\Big(\bigcap_{ p \ge M}B_H^{(p)}\,|\,\Omega_n(u,v)\Big)\ge\prod_{M\le p \le N} \Big(1-\frac{s}{p^2}\Big)-C_M \, s\, \sum_{p>N}\frac{1}{p^2} ,
\]
and conclude, upon letting $N \to \infty$, that
\[
\liminf_{n \to \infty} \P_n\Big(\bigcap_{ p \ge M}B_H^{(p)}\,|\,\Omega_n(u,v)\Big)\ge \prod_{p\ge M} \Big(1-\frac{s}{p^2}\Big) ,
\]
which finishes the proof.
\end{proof}

\subsubsection{Probability distribution of\/ $Z_M^{\star}$}\label{section:pob distr of Z and Z*}

We derive now the probability distribution of $Z_M^{\star}$ by means of  Lemma \ref{lemma:a la caibach}.

From Lemma \ref{lemma:a la caibach} and the expression \eqref{eq:indicadora en Omegan3}, we deduce, for each $(u,v)$ such that $0\le u,v <P_M$, that
\[
\begin{aligned}
\lim_{n \to \infty} \P_n(Z_M=r\, |\,\Omega_n(u,v))
&=\sum_{s=r}^{\Phi_M(u,v)} (-1)^{s-r} \binom{s}{r}  \sum_{\substack{H \subset \A_M(u,v)
\\
|H|=s}} \prod_{p \ge M}\Big(1-\frac{s}{p^2}\Big)
\\&=\sum_{s=r}^{\Phi_M(u,v)} (-1)^{s-r} \binom{s}{r} \binom{\Phi(u,v)}{s} \prod_{p \ge M} \Big(1-\frac{s}{p^2}\Big).
&\end{aligned}
\]

Therefore, from total probability and \eqref{eq:lim prob Omega_n},  we finally conclude  that
\[
\begin{aligned}
\P(Z^{\star}_M=r)&:=
\lim_{n \to \infty} \P_n(Z_M=r)
\\
&=
\frac{1}{P_M^2} \sum_{0\le u,v<P_M} \sum_{s=r}^{\Phi_M(u,v)} (-1)^{s-r} \binom{s}{r}\binom{\Phi_M(u,v)}{s} \prod_{p \ge M} \Big(1-\frac{s}{p^2}\Big)
\\
&=
\frac{1}{P_M^2} \sum_{0\le u,v<P_M} \sum_{s=r}^{M^2} (-1)^{s-r} \binom{s}{r}\binom{\Phi_M(u,v)}{s} \prod_{p \ge M} \Big(1-\frac{s}{p^2}\Big).
\end{aligned}
\]

Part of the conclusion is that each of the limits $\lim_{n \to \infty} \P_n(Z_M=r)$  exists; the other, of course, is the precise values of those limits.
The values of these limits define the probability distribution of a variable~$Z_M^\star$ with values in $\{0,\ldots, M^2\}$.

The sum of the above probabilities as $r$ runs from $r=0$ to $r=M^2$ is 1, as it should, simply because $\sum_{r=0}^s (-1)^{s - r} \binom{s}{r}=0$, unless $s=0$, in which case it is 1.

\begin{thm}
\label{teor:distribution of ZM}
With the notations above, including that recorded in \eqref{eq:formula of gammas}, we have for $M \ge 1$ and $0 \le r \le M^2$ that
\begin{equation}\label{eq:prob distr de ZMstar}
\P(Z^{\star}_M=r)=\sum_{s=r}^{M^2} (-1)^{s-r} \binom{s}{r}\,\xi(M,s) \prod_{p \ge M} \Big(1-\frac{s}{p^2}\Big).
\end{equation}
\end{thm}

Thus, as $n \to \infty$, the variable $Z_M$ in the probability space $(\N^2_n, \P_n)$  converges in distribution to the random  variable $Z_M^\star$ with probability mass function given by \eqref{eq:prob distr de ZMstar}. As we have already pointed out, the case $M = 1$ is Dirichlet's density theorem: the variable $Z^\star_1$ is a Bernoulli variable with success parameter $\P(Z_1^\star=1)=\prod_{p} (1-{s}/{p^2})=1/\zeta(2)$.

\medskip

The formula in \eqref{eq:prob distr de ZMstar} give the following probabilities, rounded to two decimal places,  for the case $M=2$:
\[\text{\footnotesize$
\begin{array}{c|cccccccccccccccccccc}
r&0&1&2&3&4
\\ \hline
\P(Z^{\star}_2=r)&0.21\%& 6.59\% & 43.00\%& 50.20\%& -
\end{array}$}
\]
The values for the case $M=3$ are:
\[\text{\footnotesize$
\begin{array}{c|cccccccccccccccccccc}
r&0&1&2&3&4&5&6&7&8&9&
\\ \hline
\P(Z^{\star}_3=r)&0.00\%& 0.02\% & 0.48\%& 4.74\%& 16.21\%& 24.41\%& 35.20\%& 17.71\%& 1.23\%& -
\end{array}$}
\]
And for the case $M=4$,
\[\text{\footnotesize$
\begin{array}{c|cccccccccccccccccccc}
r&0&1&2&3&4&5&6&7&8&
\\ \hline
\P(Z^{\star}_4=r)&0.00\%& 0.00\% & 0.00\%& 0.00\%& 0.00\%& 0.01\%\%& 0.27\%& 2.37\%& 10.67\%
\\
r&9&10&11&12&13&14&15&16
\\ \hline
\P(Z^{\star}_4=r)&25.83\%& 35.68\% & 22.18\%& 2.99\%& -& -& -& -
\end{array}$}
\]

In all these tables, a dash ``$-$'' means probability (exactly) 0: not attainable values. \label{pag:attainable}
 These are the graphical representations of the mass functions for $M=3,4,5$ in a common range $\{0,1,\dots,25\}$:
\begin{center}
\begin{tabular}{ccc}
\resizebox{4.6cm}{!}{\includegraphics{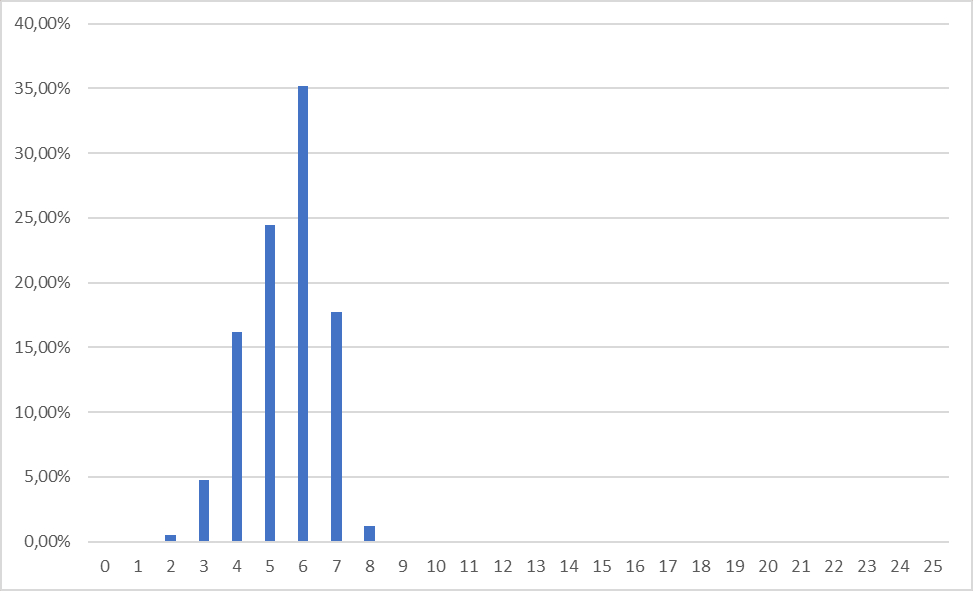}}&
\resizebox{4.6cm}{!}{\includegraphics{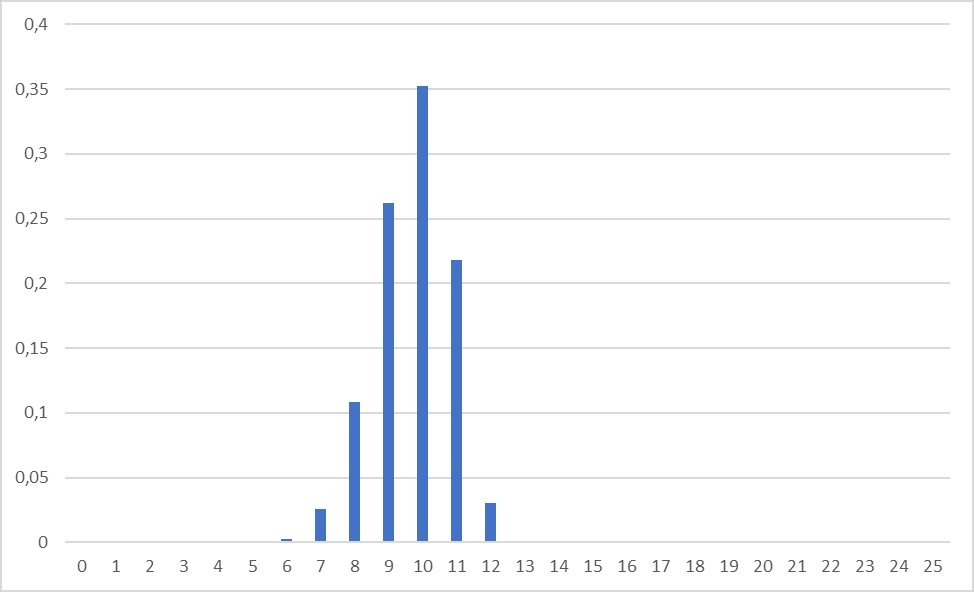}}&
\resizebox{4.6cm}{!}{\includegraphics{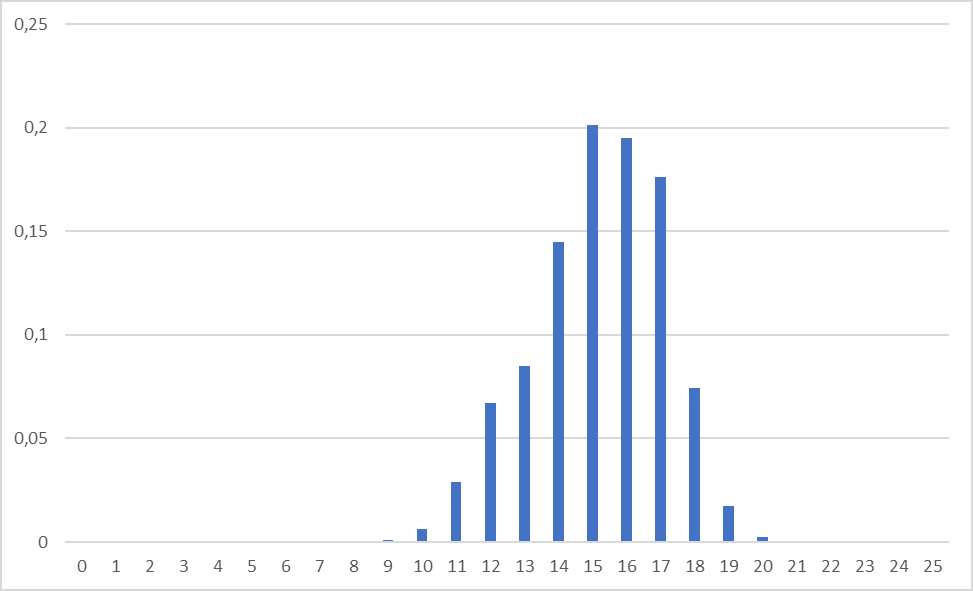}}
\\
\text{\footnotesize Mass function of $Z^*_3$.}
&
\text{\footnotesize Mass function of $Z^*_4$.}
&
\text{\footnotesize Mass function of $Z^*_5$.}
\end{tabular}
\end{center}

For the mean of $Z_M^\star$ we have, using \eqref{eq:lim mean of ZM in Nn}, that
\begin{equation}\label{eq:mean of ZMstar}
\begin{aligned}\E(Z_M^\star)&=\sum_{r=0}^{M^2} r \,\P(Z_M^\star=r)
=\sum_{r=0}^{M^2} r \lim_{n \to \infty}\P_n(Z_M=r)
\\
&=\lim_{n \to \infty}\sum_{r=0}^{M^2} r\,\P_n(Z_M=r)=\lim_{n \to \infty} \E_n(Z_M)=\frac{M^2}{\zeta(2)} \cdot
\end{aligned}
\end{equation}

For the probability generating function $G_{Z_m^\star}$ of $Z^\star_M$, we obtain immediately from \eqref{eq:prob distr de ZMstar} and the binomial theorem the following. See also Lemma~\ref{lemma:inclusion/exclusion and pgf}.
\begin{cor}\label{cor:pgf of ZM}For $|z|\le 1$,
\begin{equation}\label{eq:pgf of ZMstar}
G_{Z_m^\star}(z):= \sum_{r=0}^{M^2} \P(Z^\star_M=r) \,z^r=\sum_{s=0}^{M^2} (z-1)^s\,
\xi(M,s) \,\prod_{p\ge M}\Big(1-\frac{s}{p^2}\Big).
\end{equation}
\end{cor}

For $z=1$, both sides of the expression  \eqref{eq:pgf of ZMstar} of the probability generating  function of~$Z_M^\star$ give 1.
Differentiating $G_{Z_m^\star}(z)$  
 and evaluating 
  at $z=1$, we get 
\[
\E(Z^\star_M)=\xi(M,1) \prod_{p\ge M}\Big(1-\frac{1}{p^2}\Big),
\]
and so,  \eqref{eq:mean of ZMstar} yields
\begin{equation}
\label{eq:formula gamma(M,1)}\xi(M,1)= M^2 \prod_{p <M} 
 \Big(1-\frac{1}{p^2}\Big) .
\end{equation}

In general, by repeated differentiation of $G_{Z_m^\star}(z)$, we get, for the factorial moments,
\begin{align*}
\E\Big(\binom{Z_M^\star}{s}\Big)&=
\E\Big(\frac{1}{s!}\,(Z_M^\star(Z_M^\star-1)\cdots (Z_M^\star-s+1)\Big)&=\xi(M,s) \prod_{p \ge M} \Big(1-\frac{s}{p^2}\Big)
\\&=\frac{M^2}{\zeta(2)}\,\frac{\xi(M,s)}{\xi(M,1)}\quad\text{for $0 \le s \le M^2$.}
\end{align*}

\begin{remark}\label{remark:ZMstar is poisson} For each 
 $M$, it appears that $\Phi_M(a,b)$ takes few different values concentrated around its mean value $\xi(M,1)$ given by \eqref{eq:formula gamma(M,1)}.

For example, we have that $\xi(3,1)=6\mbox{.}8$, while the values of $\Phi_3(a,b)$ range from 5 to 8; for $M=4$, $\xi(4,1)=10\mbox{.}7$, and $\Phi_4(a,b)$ takes values between 9 and 12. See Figures~\ref{fig:cuadrados3} and~\ref{fig:cuadrados4}.

Assume \emph{for the sake of the argument} that $\Phi_M(a,b)$ \emph{were to be constantly} the integer $\xi(M,1)$.
If this were the case, then we would have, from \eqref{eq:formula of gammas},
\[
\xi(M,s)=\binom{\xi(M,1)}{s} ,\quad \mbox{for each $s\ge 0$} ,
\]
and thus, plugging this into \eqref{eq:prob distr de ZMstar}, the probability $\P(Z_M^\star=r)$ would be, approximately,
\[
\sum_{s=r}^{\xi(M,1)} (-1)^{s - r} \binom{s}{r} \binom{\xi(M,1)}{s} \prod_{p\ge M} \Big(1-\frac{s}{p^2}\Big) , \quad \mbox{for $0\le r \le M^2$}.
\]
Further, approximate
\[
\prod_{p\ge M} \Big(1-\frac{s}{p^2}\Big) \approx  \prod_{p\ge M} \Big(1-\frac{1}{p^2}\Big)^s , \quad \mbox{for each $s \ge 0$} ,
\]
and denote
\[
q_M:=\prod_{p\ge M} \Big(1-\frac{1}{p^2}\Big).
\]
With all this, 
 for $0\le r \le M^2$, $\P(Z_M^\star=r)$ would be, approximately,
\[
 \sum_{s=r}^{\xi(M,1)} (-1)^{s - r} \binom{s}{r} \binom{\xi(M,1)}{s} \,q_M^s,
\]
and then $Z_m^\star$ would follow a binomial distribution with $\xi(M,1)$ repetitions and success probability  $q_M$.

Or better, since $\xi(M,1) q_M={M^2}/{\zeta(2)}$ (see \eqref{eq:formula gamma(M,1)}), we could expect 
 the variable $Z_M^\star$ to be relatively close to a Poisson distribution with parameter $M^2/\zeta(2)$, with closeness depending  upon the factorization of $M$.

\end{remark}

\section{Coprimality correlation structure}\label{sec:correlation}

This section is devoted to study  the correlation structure between 
 the random variables that register coprimality of pairs of points in $\mathbb{N}^2$. This will be used later to obtain estimates of the variance of $Z_M^\star$.


\begin{lemma}\label{lemma:correlation}
For $(a,b) \in \N^2$ and integers $i,j \ge 0$, we have that
\[
V(a,b)\,V(a+i,b+j)=\prod_{p\,\mid \,\gcd(i,j)} (1-I_p(a)\,I_p(b)) \prod_{p \,\nmid \,\gcd(i,j)}(1-I_p(a)\,I_p(b)-I_p(a+i)\,I_p(b+j)).
\]
\end{lemma}
\begin{proof}
It follows from writing
\[
\begin{aligned}
V(a,b)&\, V(a+i,b+j)
=\prod_{p} (1-I_p(a)\,I_p(b)) \,(1-I_p(a+i)\,I_p(b+j))
\\
&\qquad=\prod_{p} (1-I_p(a)\,I_p(b)-I_p(a+i)\,I_p(b+j)+I_p(a)\,I_p(b) \,I_p(a+i)\,I_p(b+j)) ,
\end{aligned}
\]
and from  observing that if $p \mid i$ then $I_p(a)=I_p(a+i)$, and that if $p \nmid i$ then $I_p(a)\,I_p(a+i)=0$, and analogously for $b$ and $j$.
\end{proof}

For $i,j \ge 0$ fixed, consider $(a,b) \in \N_n^2 \mapsto V(a,b)V(a+i,b+j)$; it follows from Lemma~\ref{lemma:correlation}, that
\[
\lim_{n \to \infty} \E_n\big(V(\cdot, \cdot) \,V(\cdot+i,\cdot+j)\big)=\prod_{p \,\mid \,\gcd(i,j)} \Big(1-\frac{1}{p^2}\Big)\, \prod_{p \,\nmid \,\gcd(i,j)} \Big(1-\frac{2}{p^2}\Big)=\mathbf{F}\,  \Upsilon (\gcd(i,j)).
\]

In general, and analogously, we have the following.
\begin{prop}
For $i,j \ge 0$ and $k,l \ge 0$, we have  that
\begin{equation}\label{eq:precovariance}
\lim_{n \to \infty} \E_n\big(V(\cdot+k, \cdot+l) \,V(\cdot+i,\cdot+j)\big)=\mathbf{F} \, \Upsilon (\gcd(|i-k|,|j-l|)).
\end{equation}
Consequently, 
for the coefficient of correlation, 
we find that
\begin{align}\nonumber
\lim_{n \to \infty} &
\rho_n (V(\cdot+k, \cdot+l),  V(\cdot+i,\cdot+j))
\\ \label{eq:formula correlacion}
&\qquad =\frac{\zeta(2)^2\, \mathbf{F}\, \Upsilon (\gcd(|i-k|,|j-l|)) -1 }{\zeta(2)-1}=:\rho(i,j,k,l).
\end{align}
\end{prop}

Recall that $\gcd(0,0)=0$ and that $\Upsilon(0)=1/(\mathbf{F} \zeta(2)$, and observe that \eqref{eq:precovariance} 
 gives that
\[\lim_{n \to \infty} \E_n\big(V(\cdot+k, \cdot+l) V(\cdot+k,\cdot+l)\big)=\mathbf{F} \Upsilon (0)=\frac{1}{\zeta(2)} ,\]
as we already know, see \eqref{eq:dirichlet general kl}.

On account of the bounds for the function $\Upsilon$ in \eqref{eq:cotas de Upsilon}, we see that all the 
 (limit) coefficients of correlation $\rho(i,j,k,l)$ defined in \eqref{eq:formula correlacion} 
 satisfy
\[
\frac{\zeta(2)^2\,\mathbf{F}-1}{\zeta(2)-1}
 \le \rho(i,j,k,l)\le 1.\]
This (somewhat mysterious) lower bound, with numerical value $\approx -0{.}19694$, is attained when $\Upsilon(1)=1$ is plugged into \eqref{eq:formula correlacion}, that is, when $\gcd(|i-k|,|j-l|)=1$; and this happens with probability $1/\zeta(2)\approx 60\mbox{.}79\%$. The second most likely value is $\approx 0\mbox{.}4799$, and
  occurs when~$\Upsilon$ takes the value $3/2$. This happens (see \eqref{es:def de Upsilon}, and also the list of values
  in~\eqref{es:values of Upsilon}) when $\gcd(|i-k|,|j-l|)$ is   a power of 2. As the probability that $\gcd(a,b)=k$, for integers~$a$ and~$b$, is $1/(k^2\zeta(2))$, a quick calculation gives that the value $\approx 0\mbox{.}4799$ is taken with probability $1/(3\zeta(2))\approx 20\mbox{.}26\%$.

In fact, as we verify next, the (asymptotic) \emph{average correlation} is 0.
\begin{prop}
\begin{equation}\label{eq:average of correlations}
\lim_{N\to\infty} \frac{1}{N^4} \sum_{\substack{(k,i), (l,j)\in \K_N}}\rho(i,j,k,l) = 0. 
\end{equation}
\end{prop}
\begin{proof}
For integer $N$, consider
\begin{equation}\label{eq:def of AN}
A_N:=\sum_{\substack{(k,i), (l,j)\in \K_N}} \Upsilon(\gcd(|i-k|, |j-l|)).
\end{equation}
Classify now $(k,i) \in \mathcal{K}_{N}$ according to whether $k<i$, $k=i$, or 
 $k >i$, and the same with $(l,j)$, to 
  obtain that
\[
A_N=\frac{N^2}{\zeta(2)}
+4 \,\Big(N\sum_{c=1}^{N}  (N-c)\,\Upsilon(c) +\sum_{1 \le c,d \le N} (N-c)(N-d)\, \Upsilon (\gcd(c,d))\Big).
\]
Recall, \eqref{eq:cotas de Upsilon}, that the function $\Upsilon $ satisfies $1\le \Upsilon (n)\le \frac{1}{\zeta(2) \mathbf{F}}$, for any $n \ge 1$. Therefore,  we have that
\[
N \sum_{c=1}^{N}  (N-c)\,\Upsilon(c)\le \frac{1}{\zeta(2)\, \mathbf{F}} \,N^3 ,
\]
and thus,
\begin{equation}\label{eq:formula expect of ZMstar cuad}
A_N=
4 \sum_{1 \le c,d \le N} (N-c)(N-d) \,\Upsilon (\gcd(c,d))+O(N^3).
\end{equation}
As shown in 
 \eqref{eq:promedio (1-x)(1-y) upsilon},  
 we have that
\[
\frac{1}{N^2}
\sum_{1 \le c,d \le N} \Big(1-\frac{c}{N}\Big)\Big(1-\frac{d}{N}\Big) \,\Upsilon (\gcd(c,d))
\to \frac{1}{4}\,\frac{1}{\zeta(2)^2\,\mathbf{F}}, \quad\text{as $N \to \infty$}.
\]
and thus, from \eqref{eq:formula expect of ZMstar cuad}, 
\begin{equation}\label{eq:limit AN}
\frac{A_N}{N^4}\to \frac{1}{\zeta(2)^2\,\mathbf{F}}, \quad\text{as $N \to \infty$}.
\end{equation}
Recalling the definition in \eqref{eq:formula correlacion}, this gives \eqref{eq:average of correlations},
as announced.
\end{proof}


\section{Behaviour 
 as the size of the window tends to \texorpdfstring{$\infty$}{infinity}}\label{sec:limits}

\subsection{Asymptotic behaviour of the variance}

We already know, see \eqref{eq:mean of ZMstar},  that
\[
\E(Z^\star_M)^{2}=\frac{M^4}{\zeta(2)^{2}}\cdot
\]
We shall now use the results from the previous section 
to obtain that the  variance of $Z^\star_M$ is $o(M^4)$,
 as $M \to \infty$. 
\begin{prop}\label{prop:OH bound for variance} With the usual notations,
\begin{equation}\label{eq:big Oh varianza}\lim_{M\to \infty}\frac{\E((Z^\star_M)^2)}{\E(Z^\star_M)^2}=1.\end{equation}
and so,
\begin{equation}\label{eq:big Oh varianza bis}\lim_{M\to \infty}\frac{\V(Z^\star_M)}{\E(Z^\star_M)^2}=0.\end{equation}\end{prop}

\begin{proof}
 Since \[Z_M(a,b)=\sum_{(k,l) \in \K_M} V(a+k,b+l),\] we have that
\[
Z_M(a,b)^2= \sum_{\substack{(k,l), (i,j)\in \K_M}} V(a+k,b+l)\,V(a+i, b+j) ,
\]
and thus, arguing as in \eqref{eq:mean of ZMstar}, using \eqref{eq:precovariance} and that $(k,l), (i,j)\in \K_M$ simply means that $1\le k,l,i,j\le M$, we deduce that
\begin{align*}
\E((Z_M^\star)^2)&=\lim_{n \to \infty} \E_n(Z_M^2)
=\sum_{\substack{(k,l), (i,j)\in \K_M}} \mathbf{F} \,\Upsilon(\gcd(|i-k|, |j-l|))
\\
&=\mathbf{F} \sum_{\substack{(k,i), (l,j)\in \K_M}} \Upsilon(\gcd(|i-k|, |j-l|)) =\mathbf{F}\, A_M,
\end{align*}
using the notation of \eqref{eq:def of AN} in the last equality.
%
%
%
%
Finally, thanks to \eqref{eq:limit AN},
\[
\frac{\E((Z^\star_M)^2)}{\E(Z^\star_M)^2}=\frac{\zeta(2)^2}{M^4}\, \mathbf{F} A_M \to 1 , \quad \mbox{as $M \to \infty$}.\qedhere
\]
%
%
\end{proof}

\subsection{Asymptotic distribution}
For $M \ge 1$, let $Y_M^\star$ be the variable  \[Y_M^\star=\frac{Z^\star_M}{M^2}\] which  registers \emph{the average number} of coprime pairs  in a random window of side length~$M$.

We have that $\E(Y^\star_M)=1/\zeta(2)$, see once more \eqref{eq:mean of ZMstar}, and, because of
 Proposition \ref{eq:big Oh varianza}, that $\lim_{M\to \infty}\V(Y^\star_M)=0$.
Chebyshev's inequality gives immediately the following.
\begin{thm}\label{teor:limit Y_n/n^2}The random variable $Y^\star_M$ tends, {\upshape in probability}, to the constant $1/\zeta(2)$ as $M\to \infty$.
\end{thm}

One could expect a result of asymptotic normality, as $M\to\infty$, for a (convenient) normalization of the variable $Y_M^*$.
 For instance, one could consider the  variables $U_M^\star$, $M \ge 1$, given by
\[
U_M^\star=M\Big(Y_M^\star-\frac{1}{\zeta(2)}\Big).
\]
(Observe that this normalization suggests that the variance of $Z_N^*$ is of the order of $M^2$.)

But, as observed by Sugita and Takanobu in \cite{SugitaTakanobu}, and numerical experiments readily confirm, the behaviour of $U_M^\star$ may depend of the arithmetical properties of $M$. In Theorem~6 of~\cite{SugitaTakanobu}, Sugita and Takanobu obtain a description of the limit points (not a unique one) of the sequence $(U_M^\star)_{M\ge 1}$ in some $L^2$ space of the adelic framework.

\subsection{Further questions}
Here are a few questions to understand further the peculiar dependence upon $M$ of the distribution of $Z_M^\star$ or of $U_M^\star$.


\vskip1pt

(1) Is it the case that $\V(Z_M^\star)=O(M^2)$, with an absolute $O$,  improving Proposition~\ref{prop:OH bound for variance}, and as suggested by Theorem 6 in~\cite{SugitaTakanobu}?

\vskip1pt

(2) Are  the normalized variables $U_M^\star$ approximately a standard normal variable, for an appropriate sequence of sizes $M$ tending to $\infty$?

\vskip1pt

(3) Recall Remark~\ref{remark:ZMstar is poisson}. Does the total variation distance between $Z_M^\star$ and a Poisson variable with parameter $M^2/\zeta(2)$, depend upon the prime factorization of $M$?

\vskip1pt

(4)
Recall, from Remark \ref{remark:caso Z=0}, that for each $M$, there are $M$-windows which  are fully ``invisible'': all points of the window have coordinates that are not coprime, that is, $Z_M(a,b)=0$.
 These invisible windows are rare, though. In the same vein as the previous question, one would expect $\P(Z_M=0)$ to be comparable to $e^{-M^2/\zeta(2)}$.


\vskip1pt

(5) A number of possible and natural extensions of the results of this paper could
 be explored. For instance, to higher dimensions, where one would have to distinguish between fully coprime tuples and pairwise coprime tuples (and also intermediate notions of coprimality, see Section 4 of \cite{FFracsam}, or   \cite{FFcodivisibility}). And instead of the proportion  of coprime pairs, one could consider the average gcd of the pairs in the random square or other moments of gcds within the square.


\begin{thebibliography}{99}
\bibitem{Apostol}
    Apostol, T.\,M.:
    \emph{\href{https://doi.org/10.1007/978-1-4757-5579-4}{Introduction to analytic number theory}}.
    Undergrad. Texts Math., Springer, New York-Heidelberg, 1976. 

\bibitem{CaiBach}
    Cai, J.-Y. and Bach, E.:
    \href{https://doi.org/10.1007/3-540-44679-6\_53}{On testing for zero polynomials by a set of points with bounded precision}.
    In \emph{Computing and combinatorics COCOON 2001}, pp. 473--482.
    Lecture Notes in Comput. Sci. 2108, Springer, Berlin, 2001.  

\bibitem{Cesaro}
    \textsc{Ces\`{a}ro, E.}:
    \href{https://doi.org/10.1007/bf02420800}{\'{E}tude moyenne du plus grand commun diviseur de deux nombres}.
    \emph{Ann. Mat. Pura Appl.~(2)} \textbf{13} (1885), 233--268.

\bibitem{CFF}
    Cilleruelo, J., Fern\'andez, J.\,L. and Fern\'andez, P.:
    \href{https://doi.org/10.1016/j.ejc.2018.08.004}{Visible lattice points in random walks}.
    \emph{European J. Combin.} \textbf{75} (2019), 92--112. 

\bibitem{Dirichlet}
    \textsc{Dirichlet, P.\,G.\,L.}:
    \href{https://doi.org/10.1017/cbo9781139237345.007}{\"Uber die Bestimmung der mittleren Werthe in der Zahlentheorie}.
    In \emph{Abhandlungen der K\"oniglich Preussischen Akademie der Wissenschaften con 1849}, 69--83.

\bibitem{Feller}
    Feller, W.:
    \emph{\href{https://www.wiley.com/en-us/An+Introduction+to+Probability+Theory+and+Its+Applications,+Volume+1,+3rd+Edition-p-9780471257080}{An introduction to probability theory and its applications. {V}ol. I}}.
    Third edition.
    John Wiley \& Sons, New York-London-Sydney, 1968.  

\bibitem{FFequi1}
    Fern\'andez, J.\,L. and Fern\'andez, P.:
    Equidistribution and coprimality.
    Preprint 2013, arXiv:\,\href{https://arxiv.org/abs/1310.3802}{1310.3802}.

\bibitem{FFcodivisibility}
    Fern\'andez, J.\,L. and Fern\'andez, P.:
    {R}andom index of codivisibility.
    Preprint 2013, arXiv:\,\href{https://arxiv.org/abs/1310.4681}{1310.4681}.

\bibitem{FFracsam}
    Fern\'andez, J.\,L. and Fern\'andez, P.:
    \href{https://doi.org/10.1007/s13398-020-00960-x}{Divisibility properties of random samples of integers}.
    \emph{Rev. R. Acad. Cienc. Exactas F{\'{\i}}s. Nat. Ser. A Mat. RACSAM} \textbf{115} (2021), no.~1, article no.~26, 35 pp. 

\bibitem{FF}
    Fern{\' a}ndez, J.\,L. and Fern{\' a}ndez, P.:
    \href{https://doi.org/10.37236/11424}{Some arithmetic properties of P\'olya's urn}.
    \emph{Electron. J. Comb.} \textbf{30} (2023), no. 2, article no. P2.11, 32 pp. 

\bibitem{Gerber}
    Gerber, H.\,U.:
    A proof of the Schuette--Nesbitt formula for dependent events.
    \emph{Actuarial Research Clearing House} \textbf{1} (1979), 9--10.

\bibitem{Gerberbook}
    Gerber, H.\,U.:
    \emph{Life insurance mathematics}. Third edition, Springer, Berlin, 1997.

\bibitem{HardyWright}
    Hardy, G.\,H. and Wright, E.\,M.:
    \emph{\href{https://doi.org/10.1093/oso/9780199219858.001.0001}{An introduction to the theory of numbers}}.
    Sixth edition.
    Oxford University Press, Oxford, 2008. 

\bibitem{HerzogStewart}
    Herzog, F. and Stewart, B.\,M.:
    \href{https://doi.org/10.2307/2317753}{Patterns of visible and nonvisible lattice points}.
    \emph{Amer. Math. Monthly} \textbf{78} (1971), 487--496. 

\bibitem{Martineau}
    Martineau, S.:
    \href{https://doi.org/10.1214/21-ecp381}{On coprime percolation, the visibility graphon, and the local limit of the {GCD} profile}.
    \emph{Electron. Commun. Probab.} \textbf{27} (2022), article no.~8,  14 pp. 

\bibitem{SugitaTakanobu}
    Sugita, H. and Takanobu, S.:
    \href{https://projecteuclid.org/journals/osaka-journal-of-mathematics/volume-40/issue-4/The-probability-of-two-integers-to-be-co-prime-revisited/ojm/1153493406.full}
    {The probability of two integers to be co-prime, revisited --\,on the behavior of {CLT}-scaling limit}.
    \emph{Osaka J. Math.} \textbf{40} (2003), no.~4, 945--976. 

\end{thebibliography}
\end{document}